\documentclass[12pt,leqno]{amsart}
\usepackage{amsmath,amsthm,amsfonts,amsbsy,multicol,fancybox}
\usepackage[utf8]{inputenc}
\usepackage[T1]{fontenc}

\usepackage{graphicx}\usepackage{amssymb,latexsym}
\usepackage{epstopdf}
\usepackage{color}
\DeclareGraphicsExtensions{.pdf, .jpg, .pnp, .bmp, .fig}
\usepackage{tikz}
\usepackage[all]{xy}
\input xy
\xyoption{all}

\usepackage[hidelinks]{hyperref}
\usepackage{amssymb,latexsym}

\usepackage{thmtools}

\newcounter{argument}
\newenvironment{argument}[1][\medskip]{%
\refstepcounter{argument}
\par\medskip
\noindent\phantomsection
\textbf{#1~\thesection.\arabic{argument}\,\,}\rmfamily\em}{\hspace{\fill}$\Box$\par\smallskip\noindent}

\newcommand{\bass}{\begin{argument}[Assumption]}\newcommand{\eass}{\end{argument}}
\newcommand{\bth}{\begin{argument}[Theorem]} \newcommand{\ethe}{\end{argument}}
\newcommand{\bre}{\begin{argument}[Remark]}      \newcommand{\ere}{\end{argument}}
\newcommand{\ble}{\begin{argument}[Lemma]}       \newcommand{\ele}{\end{argument}}
\newcommand{\bde}{\begin{argument}[Definition]}   \newcommand{\ede}{\end{argument}}
\newcommand{\bco}{\begin{argument}[Corollary]}     \newcommand{\eco}{\end{argument}}
\newcommand{\bpr}{\begin{argument}[Proposition]}  \newcommand{\epr}{\end{argument}}
\newcommand{\bexam}{\begin{argument}[Example]}\newcommand{\eexam}{\end{argument}}

\newcommand{\bpf}{\begin{proof}}\newcommand{\epf}{\end{proof}}
\newcommand{\barr}{\begin{array}}\newcommand{\earr}{\end{array}}
\newcommand{\beao}{\begin{eqnarray*}}\newcommand{\eeao}{\end{eqnarray*}\noindent}
\newcommand{\beam}{\begin{eqnarray}}\newcommand{\eeam}{\end{eqnarray}\noindent}
\newcommand{\beqq}{\begin{equation}}\newcommand{\eeqq}{\end{equation}\noindent}

\newcommand{\ov}{\overline} \newcommand{\un}{\underbrace}
\newcommand{\wt}{\widetilde}
\newcommand{\wh}{\widehat}

\newcommand{\D}{\Delta}
\newcommand{\vep}{\varepsilon}  \newcommand{\ep}{\epsilon}

\newcommand{\w}{\omega} \newcommand{\W}{\Omega}

\newcommand{\bbc}{{\mathcal C}} 

\newcommand{\bfE}{{\mathbb E}}\newcommand{\bbE}{{\mathcal E}} 
\newcommand{\bbf}{{\mathcal F}}

\newcommand{\bbi}{{\mathbb I}}\newcommand{\bbI}{{\mathcal I}}

\newcommand{\bbl}{{\mathcal L}}
\newcommand{\bbm}{{\mathcal M}}
 \newcommand{\bbN}{{\mathbb N}}
\newcommand{\bbo}{{\mathcal O}}
\newcommand{\bfP}{{\mathbb P}}
                              
 \newcommand{\bbR}{{\mathbb R}}

\graphicspath{{./../Images/}}

\begin{document}

\title[A boundary preserving numerical scheme for the Wright-Fisher model]{A boundary preserving numerical scheme for the Wright-Fisher model}

\author[I. S. Stamatiou]{I. S. Stamatiou}
\email{joniou@gmail.com}

\begin{abstract}
We are interested in the numerical approximation of non-linear stochastic differential equations (SDEs) with solution in a certain domain. Our goal is to construct explicit numerical schemes that preserve that structure. We generalize the semi-discrete method \emph{Halidias N. and Stamatiou I.S. (2016), On the numerical solution of some non-linear stochastic differential equations using the  Semi-Discrete method, Computational Methods in Applied Mathematics,16(1)}  and propose a numerical scheme, for which we prove a strong convergence result, to a class of SDEs that appears in population dynamics and ion channel dynamics within cardiac and neuronal cells. We furthermore extend our scheme to a multidimensional case.   
\end{abstract}

\date\today

\maketitle

\section{Introduction}\label{NSF:sec:intro}
\setcounter{equation}{0}

Let  $T>0$ and $(\Omega, \bbf, \{\bbf_t\}_{0\leq t\leq T}, \bfP)$ be a complete probability space and let $W_{t,\w}:[0,T]\times\W\rightarrow\bbR$ be a one-dimensional Wiener process adapted to the filtration  $\{\bbf_t\}_{0\leq t\leq  T}.$
  We are interested in the numerical approximation of the following scalar stochastic differential equation (SDE),
\beqq  \label{NSF-eq:WFmodel}
x_t =x_0 + \int_0^t (k_1 - k_2x_s)ds + k_3\int_0^t \sqrt{x_s(1-x_s)}dW_s,
\eeqq
where $k_i>0, i=1,2,3.$ A boundary classification result, see Appendix \ref{NSF-ap:Boundary Classification}, 
implies that $0<x_t<1$ a.s. when $x_0\in(0,1)$ and $0<k_1<k_2.$  We therefore aim for a numerical scheme which apart from strongly converging to the true solution of (\ref{NSF-eq:WFmodel}), produces values in the same domain, i.e. in  $(0,1).$ In other words, we are interested in numerical schemes that have an \emph{eternal life time}.

\bde[Eternal Life time of numerical solution]
Let $D\subseteq\bbR^d$ and consider a process $(X_t)$ well defined on the domain $\ov{D},$ with initial condition $X_0\in\ov{D}$ and such that 
$$\bfP(\{\w\in\W: X(t,\w)\in\ov{D}\})=1,$$ for all $t>0.$ A numerical solution $(Y_{t_n})_{n\in\bbN}$ has an \emph{eternal life time} if 
$$\bfP(Y_{n+1}\in \ov{D}\, \big| \,Y_{n}\in \ov{D})=1. $$
\ede

In \cite{schurz:1996} the main interest is in the domain $D=\bbR_+.$  Moreover, it is clear that the Euler-Maruyama scheme has always a finite life time.

The proposed semi-discrete (SD) iterative scheme for the numerical approximation of (\ref{NSF-eq:WFmodel}) reads
\beqq\label{NSF-eq:SD scheme}
y_{t_{n+1}}^{SD}=\sin^2 \left(\frac{k_3}{2}\D W_n + \arcsin(\sqrt{y_n})\right),
\eeqq
where 
$$
y_n:=y_{t_n} + \left(k_1 - \frac{(k_3)^2}{4} + y_{t_n}\left(\frac{(k_3)^2}{2}-k_2\right)\right)\cdot\D
$$
and
$\D W_n:= W_{t_{n+1}}-W_{t_n},$ are the increments of the Wiener process and the discretization step $\D$ is such that (\ref{NSF-eq:SD scheme})  is well-defined. By construction, the SD scheme (\ref{NSF-eq:SD scheme}) possesses an eternal life time. To get (\ref{NSF-eq:SD scheme}) we use an additive semi-discretization of the drift coefficient.. Briefly saying, a part of the SDE  is discretized in a certain way such that the resulting SDE to be solved has an analytical solution (see details in Section \ref{NSF-sec:setting}). This is also a special feature of the method, since in the derivation of it, instead of an algebraic equation a new SDE has to be solved. The SD method can also reproduce the Euler scheme. The semi-discrete method was originally proposed in \cite{halidias:2012}. 

An attempt in that direction, i.e. in constructing explicit numerical schemes with an eternal life time, has been made in \cite{halidias_stamatiou:2016} where a class of one-dimensional SDEs with non-negative solutions is treated,  which covers cases like that  of the Heston $3/2$-model, a popular model in the field of financial mathematics which is super-linear. The case of sub-linearities is also treated in \cite{halidias_stamatiou:2015} where the domain is still $\bbR_+.$

The purpose of this paper is to generalize further the method to preserve the structure of the original SDE. In the previous works, the suggested schemes  preserve positivity; all the quantities appearing belong to the field of finance and are meant to be non-negative. The application that motivated us  now, is used in population dynamics to describe fluctuations in gene frequency of reproducing individuals among finite populations \cite{ewens:2012} and in a different setting for the description of the random behavior of ion channels within cardiac and neuronal cells  (cf. \cite{dangerfield:2012}, \cite{goldwyn:2011}, \cite{dangerfield:2012b} and references therein). We are able in that case to preserve the domain of the original process. In fact, in applications in biology we have to solve systems of SDEs. The extension of the Wright-Fisher model to the multidimensional case has been proposed in \cite{griffiths:1979} and \cite{griffiths:1980}. Here, we will treat a three-state  system as in \cite[Sec. 6]{dangerfield:2012} given by the following system of SDEs
\beam 
\nonumber 
X_t^{(1)} &=& X_0^{(1)} + \int_0^t (k_1^{(1,1)}+k_1^{(1,2)}X_s^{(2)} - k_2^{(1)}X_s^{(1)})ds\\
\label{NSF-eq:WFmultimodel}&&+ k_3^{(1,1)}\int_0^t \sqrt{X_s^{(1)}X_s^{(2)}}dW_s^{(1)} + k_3^{(1,2)}\int_0^t \sqrt{X_s^{(1)}(1-X_s^{(1)}-X_s^{(2)})}dW_s^{(2)},\\
\nonumber 
X_t^{(2)} &=& X_0^{(2)} + \int_0^t (k_1^{(2,1)}+k_1^{(2,2)}X_s^{(1)} - k_2^{(2)}X_s^{(2)})ds\\
\label{NSF-eq:WFmultimodel2}&&+ k_3^{(2,1)}\int_0^t \sqrt{X_s^{(1)}X_s^{(2)}}dW_s^{(1)} + k_3^{(2,3)}\int_0^t \sqrt{X_s^{(2)}(1-X_s^{(1)}-X_s^{(2)})}dW_s^{(3)},
\eeam
where $X^{(i)}$ is the proportion of alleles or channels in state $i, i=1,2,$ and $1-X^{(1)}-X^{(2)}$ is the proportion in state $3.$

In Section \ref{NSF-sec:setting} we provide the setting and the main goal which concerns the mean-square convergence of the proposed structure-preserving SD scheme (\ref{NSF-eq:SD scheme}) for the approximation of a modification of (\ref{NSF-eq:WFmodel}) with dynamics described by $\wh{W}$ (see (\ref{NSF-eq:NewWiener})).  We also discuss the multidimensional case (\ref{NSF-eq:WFmultimodel})-(\ref{NSF-eq:WFmultimodel2}).

In Section \ref{ExSD-sec:Intro} we treat a more general class of SDEs. We further extend the analysis of the semi-discrete method introduced in \cite{halidias_stamatiou:2016}. We cover the sub-linear diffusion case and show as in \cite[Th. 2.1]{halidias_stamatiou:2016} the strong convergence of the proposed numerical scheme to the true solution.  

Section  \ref{NSF-sec:Experiments} is devoted to numerical experiments. The proofs of all the results are given in the sections to follow, that is in Sections \ref{ExSD-sec:Proof}  and \ref{NSF-sec:Proof}. 

\section{The setting and the main goal.}\label{NSF-sec:setting}
\setcounter{equation}{0}
Consider the partition $0=t_0<t_1<\ldots<t_N=T$ with uniform discretization step $\D=T/N$ and the following process
\beam\nonumber
y_t^{SD}&=& y_{t_n} + \int_{t_n}^{t_{n+1}}\left(k_1 - \frac{(k_3)^2}{4} + y_{t_n}\left(\frac{(k_3)^2}{2}-k_2\right)\right)ds  + \int_{t_n}^{t} \frac{(k_3)^2}{4}(1-2y_s)ds\\
\label{NSF-eq:SD process}&& + k_3\int_{t_n}^{t} \sqrt{y_s(1-y_s)}\,\textup{sgn}(z_s)dW_s,
\eeam
for $ t\in(t_n,t_{n+1}],$ with $y_0=x_0$ a.s. 
\beqq\label{NSF-eq:initial}
y_n:=y_{t_n} + \left(k_1 - \frac{(k_3)^2}{4} + y_{t_n}\left(\frac{(k_3)^2}{2}-k_2\right)\right)\cdot\D
\eeqq
and 
\beqq\label{NSF-eq:SDsgn term}
z_t= \sin \left(k_3\D W_n^t + 2\arcsin(\sqrt{y_n})\right), 
\eeqq
where $\D W_n^t:= W_{t}-W_{t_n}.$ Process (\ref{NSF-eq:SD process}) has jumps at nodes $t_n$ of order $\D$ and the solution in each step is given by,  see Appendix \ref{NSF-ap:Solution process}, 
\beqq\label{NSF-eq:SDsolution}
y_t^{SD}=\sin^2 \left(\frac{k_3}{2}\D W_n^t + \arcsin(\sqrt{y_n})\right),
\eeqq
 which has the pleasant feature that $y_t^{SD}\in(0,1)$ when $y_0\in(0,1).$ Process (\ref{NSF-eq:SDsolution}) is well defined when $0< y_n<1,$ i.e. when
\beqq\label{NSF-eq:arccos condition}
0< y_{t_n} + \left(k_1 - \frac{(k_3)^2}{4} + y_{t_n}\left(\frac{(k_3)^2}{2}-k_2\right)\right)\cdot\D<1.
\eeqq 
Therefore, we assume the following condition for the well-posedness of the SD scheme (\ref{NSF-eq:SDsolution}).

\bass\label{NSF:assA} 
Let the discretization step $\D$ be such that (\ref{NSF-eq:arccos condition}) holds.
\eass

\bre\label{NSF-rem:stepsize}
Note that in general the discretization step $\D$ satisfying  (\ref{NSF-eq:arccos condition}) is a r.v. depending on $\w.$ The $\w$-dependence is inherited through the increments $\D W_n(\w)$ which in turn affect the sequence $(y_{t_n})_{n\in\bbN}$.  Nevertheless under the assumptions on the parameters considered later on the step $\D$  is not a r.v. but a fixed sufficiently small number.
\ere

Now, we consider the process
\beqq\label{NSF-eq:NewWiener}
 \wh{W}_t:=\int_{0}^{t} \textup{sgn}(z_s)dW_s,
\eeqq
which is a martingale with quadratic variation $< \wh{W}_t, \wh{W}_t>=t$ and thus a standard Brownian motion w.r.t. its own filtration, justified by L\'evy's theorem \cite[Th. 3.16, p.157]{karatzas_shreve:1988} and consequently  (\ref{NSF-eq:SD process}) becomes
\beqq\label{NSF-eq:SD processNW}
y_t^{SD}= y_n + \int_{t_n}^{t} \frac{(k_3)^2}{4}(1-2y_s)ds + k_3\int_{t_n}^{t} \sqrt{y_s(1-y_s)}d\wh{W}_s,
\eeqq
Moreover, consider 
\beqq  \label{NSF-eq:WFmodel_NW}
\wh{x}_t =x_0 + \int_0^t (k_1 - k_2\wh{x}_s)ds + k_3\int_0^t \sqrt{\wh{x}_s(1-\wh{x}_s)}d\wh{W}_s.
\eeqq

The process $(x_t)$ of (\ref{NSF-eq:WFmodel}) and the process $(\wh{x}_t)$ of (\ref{NSF-eq:WFmodel_NW}) have the same distribution. 
Our main goal is to deduce an estimate of the form 
$$\lim_{\D\downarrow0}\bfE\sup_{0\leq t\leq T}|y_t^{SD} - x_t|^2=0.$$ 
In Theorem \ref{NSF-theorem:Pol rate} below, we deduce that 
$$\lim_{\D\downarrow0}\bfE\sup_{0\leq t\leq T}|y_t^{SD} - \wh{x}_t|^2=0.$$ 
By a simple application of  the triangle inequality  
 we deduce an analogous result for the unique solution of (\ref{NSF-eq:WFmodel}), i.e. $\lim_{\D\downarrow0}\bfE\sup_{0\leq t\leq T}|y_t^{SD}-x_t|^2=0.$ 
 We present in Appendix \ref{NSF-ap:uniform_moment_bound} the details. To simplify notation we write $\wh{W}, (\wh{x}_t)$ as $W, (x_t)$ respectively.
 
\bth\label{NSF-theorem:Pol rate}[Strong convergence]
Let Assumption \ref{NSF:assA} hold. Then, the semi-discrete scheme (\ref{NSF-eq:SD process}) converges strongly in the mean-square sense to the true solution of  (\ref{NSF-eq:WFmodel}), that is 
$$
\lim_{\D\downarrow0}\bfE\sup_{0\leq t\leq T}|y_{t}^{SD}-x_{t}|^2=0. 
$$
\ethe

As already noted in Remark  \ref{NSF-rem:stepsize}, in order to apply  Theorem \ref{NSF-theorem:Pol rate} we have to find a sufficiently small step-size  $\D$ such that  (\ref{NSF-eq:arccos condition}) holds, i.e. 
\beqq\label{NSF-eq:arccos condition_alt}
0< y_{t_n}(1 + \beta \D) + \alpha\D <1.
\eeqq
where 
\beqq\label{NSF-eq:additional_par}
\alpha:=k_1 - \frac{(k_3)^2}{4} ,\quad \beta:=\frac{(k_3)^2}{2}-k_2.
\eeqq
To simplify the conditions on $\alpha, \beta, \D$, when necessary, we may adopt the following procedure. We consider a perturbation of order $\D$ in the initial condition of (\ref{NSF-eq:SD processNW}), that is 
\beqq\label{NSF-eq:SD processNWalt}
\wt{y}_t^{SD}= \wt{y}_n + \int_{t_n}^{t} \frac{(k_3)^2}{4}(1-2\wt{y}_s)ds + k_3\int_{t_n}^{t} \sqrt{\wt{y}_s(1-\wt{y}_s)}d\wh{W}_s,
\eeqq
for $ t\in(t_n,t_{n+1}],$ with $\wt{y}_0=x_0$ a.s. and
\beqq\label{NSF-eq:initial_al}
\wt{y}_n:=\frac{\wt{y}_{t_n}(1 + \beta \D) + \alpha\D}{1  + (\alpha+\beta)\D}.
\eeqq
The following result is a consequence of Theorem \ref{NSF-theorem:Pol rate}. 
\bpr\label{NSF-prop:Pol rate}[Strong convergence]
Let $(k_3)^2<2k_2$ and $\D<-1/\beta,$ where $\beta$ is given by (\ref{NSF-eq:additional_par}). Then, the semi-discrete scheme (\ref{NSF-eq:SD processNWalt}) converges strongly in the mean-square sense to the true solution of  (\ref{NSF-eq:WFmodel}), 
that is 
$$
\lim_{\D\downarrow0}\bfE\sup_{0\leq t\leq T}|\wt{y}_{t}^{SD}-x_{t}|^2=0. 
$$
\epr

Now, we turn to the approximation of the solution of system  (\ref{NSF-eq:WFmultimodel})-(\ref{NSF-eq:WFmultimodel2}). This time we also discretize  the diffusion coefficient in a multiplicative way such that the resulting SDEs for each component can be solved analytically. In particular we consider the following processes   
\beam\nonumber
_{SD}Y_t^{(1)}&=& Y_{t_n}^{(1)} + \int_{t_n}^{t_{n+1}}\left(k_1^{(1,1)}+k_1^{(1,2)}Y_{t_n}^{(2)} - \frac{(k_3^{(1,1)})^2Y_{t_n}^{(2)} + (k_3^{(1,2)})^2(1-Y_{t_n}^{(1)}-Y_{t_n}^{(2)})}{4(1-Y_{t_n}^{(1)})}\right)ds\\
\nonumber&& + \int_{t_n}^{t_{n+1}}Y_{t_n}^{(1)}\left(\frac{(k_3^{(1,1)})^2Y_{t_n}^{(2)} + (k_3^{(1,2)})^2(1-Y_{t_n}^{(1)}-Y_{t_n}^{(2)})}{2(1-Y_{t_n}^{(1)})} -k_2^{(1)}\right)ds \\ 
\nonumber&& + \int_{t_n}^{t} \frac{(k_3^{(1,1)})^2Y_{t_n}^{(2)}+ (k_3^{(1,2)})^2(1-Y_{t_n}^{(1)}-Y_{t_n}^{(2)})}{4(1-Y_{t_n}^{(1)})}(1-2Y_s^{(1)})ds\\
\nonumber&& + k_3^{(1,1)}\sqrt{\frac{Y_{t_n}^{(2)}}{1-Y_{t_n}^{(1)}}}\int_{t_n}^{t} \sqrt{Y_s^{(1)}(1-Y_s^{(1)})}\,\textup{sgn}(z_s^{(1)})dW_s^{(1)}\\
\label{NSF-eq:SDmultiprocess1}&& + k_3^{(1,2)}\sqrt{\frac{1-Y_{t_n}^{(1)}-Y_{t_n}^{(2)}}{1-Y_{t_n}^{(1)}}}\int_{t_n}^{t} \sqrt{Y_s^{(1)}(1-Y_s^{(1)})}\,\textup{sgn}(z_s^{(1)})dW_s^{(2)},
\eeam 
\beam\nonumber
_{SD}Y_t^{(2)}&=& Y_{t_n}^{(2)} + \int_{t_n}^{t_{n+1}}\left(k_1^{(2,1)}+k_1^{(2,2)}Y_{t_n}^{(1)} - \frac{(k_3^{(2,1)})^2Y_{t_n}^{(1)} + (k_3^{(2,3)})^2(1-Y_{t_n}^{(1)}-Y_{t_n}^{(2)})}{4(1-Y_{t_n}^{(2)})}\right)ds\\
\nonumber&& + \int_{t_n}^{t_{n+1}}Y_{t_n}^{(2)}\left(\frac{(k_3^{(2,1)})^2Y_{t_n}^{(1)} + (k_3^{(2,3)})^2(1-Y_{t_n}^{(1)}-Y_{t_n}^{(2)})}{2(1-Y_{t_n}^{(2)})} -k_2^{(2)}\right)ds \\ 
\nonumber&& + \int_{t_n}^{t} \frac{(k_3^{(2,1)})^2Y_{t_n}^{(1)}+ (k_3^{(2,3)})^2(1-Y_{t_n}^{(1)}-Y_{t_n}^{(2)})}{4(1-Y_{t_n}^{(2)})}(1-2Y_s^{(2)})ds\\
\nonumber&& + k_3^{(2,1)}\sqrt{\frac{Y_{t_n}^{(1)}}{1-Y_{t_n}^{(2)}}}\int_{t_n}^{t} \sqrt{Y_s^{(2)}(1-Y_s^{(2)})}\,\textup{sgn}(z_s^{(2)})dW_s^{(1)}\\
\label{NSF-eq:SDmultiprocess2}&& + k_3^{(2,3)}\sqrt{\frac{1-Y_{t_n}^{(1)}-Y_{t_n}^{(2)}}{1-Y_{t_n}^{(2)}}}\int_{t_n}^{t} \sqrt{Y_s^{(2)}(1-Y_s^{(2)})}\,\textup{sgn}(z_s^{(2)})dW_s^{(3)},
\eeam
for $ t\in(t_n,t_{n+1}],$ with $Y_0^{(i)}=X_0^{(i)}, i=1,2,3$  a.s. and 
$$
z_t^{(1)}= \sin \left(k_3^{(1,1)}\sqrt{\frac{Y_{t_n}^{(2)}}{1-Y_{t_n}^{(1)}}}\D W_t^{(1)} + k_3^{(1,2)}\sqrt{\frac{1-Y_{t_n}^{(1)}-Y_{t_n}^{(2)}}{1-Y_{t_n}^{(1)}}}\D W_t^{(2)}  + 2\arcsin(\sqrt{y_n^{(1)}})\right), 
$$
$$
z_t^{(2)}= \sin \left(k_3^{(2,1)}\sqrt{\frac{Y_{t_n}^{(1)}}{1-Y_{t_n}^{(2)}}}\D W_t^{(1)} + k_3^{(2,3)}\sqrt{\frac{1-Y_{t_n}^{(1)}-Y_{t_n}^{(2)}}{1-Y_{t_n}^{(2)}}}\D W_t^{(3)} + 2\arcsin(\sqrt{y_n^{(2)}})\right), 
$$
where $y_n^{(i)}, i=1,2,$ are the deterministic parts of (\ref{NSF-eq:SDmultiprocess1}) and  (\ref{NSF-eq:SDmultiprocess2}) respectively and $\D W_t^{(i)}:= W_{t}^{(i)}-W_{t_n}^{(i)}.$ 
We set $Y_{t_n}^{(i)}=1-\vep, Y_{t_n}^{(j)}=\vep/2,$ for  $i,j=1,2$ with $i\neq j$  whenever  $Y_{t_n}^{(i)}>1-\vep$ for a tolerance $\vep>0.$ Finally $Y_{t_n}^{(3)}=1-Y_{t_n}^{(1)}-Y_{t_n}^{(2)}.$

Processes (\ref{NSF-eq:SDmultiprocess1}) and (\ref{NSF-eq:SDmultiprocess2}) have jumps at nodes $t_n$ of order $\D$ and their solution in each step is given respectively by,  see Appendix \ref{NSF-ap:Solution process}, 
\beqq\label{NSF-eq:SDmultisolution1}
_{SD}Y_t^{(1)}=\sin^2 \left(\frac{k_3^{(1,1)}}{2}\sqrt{\frac{Y_{t_n}^{(2)}}{1-Y_{t_n}^{(1)}}}\D W_t^{(1)} + \frac{k_3^{(1,2)}}{2}\sqrt{\frac{1-Y_{t_n}^{(1)}-Y_{t_n}^{(2)}}{1-Y_{t_n}^{(1)}}}\D W_t^{(2)}  + \arcsin(\sqrt{y_n^{(1)}})\right),
\eeqq
and 
\beqq\label{NSF-eq:SDmultisolution2}
_{SD}Y_t^{(2)} = 
\sin \left(\frac{k_3^{(2,1)}}{2}\sqrt{\frac{Y_{t_n}^{(1)}}{1-Y_{t_n}^{(2)}}}\D W_t^{(1)} + \frac{k_3^{(2,3)}}{2}\sqrt{\frac{1-Y_{t_n}^{(1)}-Y_{t_n}^{(2)}}{1-Y_{t_n}^{(2)}}}\D W_t^{(3)} + \arcsin(\sqrt{y_n^{(2)}})\right),
\eeqq
which has the pleasant feature that $Y_t^{SD}\in(0,1)^3$ when $Y_0\in(0,1)^{3}.$ Processes (\ref{NSF-eq:SDmultisolution1}) and (\ref{NSF-eq:SDmultisolution2}) are well defined when $0< y_n^{(i)}<1, i=1,2.$

Working as before, considering this time the processes
\beqq\label{NSF-eq:NewWienermulti}
 \wh{W}_t^{(i)}:=\int_{0}^{t} \textup{sgn}(z_s^{(i)})dW_s^{(i)},
\eeqq
we conclude to the following result.
 \bth\label{NSF-theorem:Pol ratemulti}[Strong convergence]
Let the discretization step be such that $0< y_n^{(i)}<1, i=1,2.$ Then, the semi-discrete scheme (\ref{NSF-eq:SDmultiprocess1})-(\ref{NSF-eq:SDmultiprocess2}) converges strongly in the mean-square sense to the true solution of  (\ref{NSF-eq:WFmultimodel})-(\ref{NSF-eq:WFmultimodel2}), that is 
$$
\lim_{\D\downarrow0}\bfE\sup_{0\leq t\leq T}||_{SD}Y_{t}-X_{t}||_2^2=0, 
$$
where $||x||_2=\sqrt{\sum_{i=1}^{d} x_i^2}$ for a $d$-dimensional vector $x$.
\ethe

We apply the above results in Section \ref{NSF-sec:Experiments} and prove them in  Section \ref{NSF-sec:Proof}.

\section{An extension of the semi-discrete method.}\label{ExSD-sec:Intro}
\setcounter{equation}{0}

Throughout, let $T>0$ and $(\Omega, \bbf, \{\bbf_t\}_{0\leq t\leq T}, \bfP)$ be a complete probability space, meaning that the filtration $ \{\bbf_t\}_{0\leq t\leq T} $ satisfies the usual conditions, i.e. is right continuous and $\bbf_0$ includes all $\bfP$-null sets. Let $W_{t,\w}:[0,T]\times\W\rightarrow\bbR$ be a one-dimensional Wiener process adapted to the filtration $\{\bbf_t\}_{0\leq t\leq  T}.$  Consider the following stochastic differential equation (SDE),
\beqq\label{ExSD-eq:general sde}
x_t=x_0 + \int_{0}^{t}a(s,x_s)ds + \int_{0}^{t}b(s,x_s)dW_s,\quad t\in [0,T],
\eeqq
where the coefficients $a,b: [0,T]\times \bbR\rightarrow\bbR$ are measurable functions such that (\ref{ExSD-eq:general sde}) has a unique strong solution and $x_0$ is independent of all $\{W_t\}_{0\leq t\leq T}.$  SDE (\ref{ExSD-eq:general sde}) has non-autonomous coefficients, i.e. $a(t,x), b(t,x)$ depend explicitly on $t.$

To be more precise, we assume the existence of a predictable stochastic process $x:[0,T]\times \W\rightarrow \bbR$ such that (\cite[Def. 2.1]{mao:2007}),
$$
\{a(t,x_t)\}\in\bbl^1([0,T];\bbR), \quad \{b(t,x_t)\}\in\bbl^2([0,T];\bbR)
$$
and
$$
\bfP\left[x_t=x_0 + \int_{0}^{t}a(s,x_s)ds + \int_{0}^{t}b(s,x_s)dW_s\right]=1, \quad \hbox{ for every } t\in[0,T].
$$
SDEs of the form (\ref{ExSD-eq:general sde}) have rarely explicit solutions, thus numerical approximations are necessary for simulations of the paths $x_t(\w),$ or for approximation of functionals of the form $\bfE F(x),$ where $F:\bbc([0,T],\bbR)\rightarrow\bbR.$ We are interested in strong approximations (mean-square) of (\ref{ExSD-eq:general sde}), in the case of super- or sub-linear drift and diffusion coefficients and cover cases not included in the previous work \cite{halidias_stamatiou:2016}.

The purpose of this section is to further generalize the  semi-discrete (SD) method  covering cases of sub-linear diffusion coefficients such as the Cox-Ingersoll-Ross model (CIR) or the Constant Elasticity of Variance model (CEV) (cf. \cite[(1.2) and  (1.4)]{halidias_stamatiou:2016}), where also an additive discretization is considered.   

\bass\label{ExSD-ass:A}  
Let $f_1(s,x):[0,T]\times\bbR\rightarrow\bbR$ and $f_2(s,r,x,y), g(s,r,x,y):[0,T]^2\times\bbR^2\rightarrow\bbR$ be such that $f_1(s,x) + f_2(s,s,x,x)=a(s,x), g(s,s,x,x)=b(s,x),$ where $f_1,f_2,g$ satisfy the following conditions
$$
f_1(s,x)\leq C(1 + |x|^l), \qquad \hbox{(Polynomial Growth)},
$$
for some appropriate $0<l$  (we take $0<l\leq p/2$ in Theorem \ref{ExSD-thm:strong_conv})
$$
|f_2(s_1,r_1,x_1,y_1) - f_2(s_2,r_2,x_2,y_2)|\leq C_R \Big( |s_1-s_2| + |r_1-r_2| + |x_1-x_2| + |y_1-y_2| \Big)
$$
and
$$
|g(s_1,r_1,x_1,y_1) - g(s_2,r_2,x_2,y_2)|\leq C_R \Big( |s_1-s_2| + |r_1-r_2|  + |x_1-x_2| + |y_1-y_2| + |x_1-x_2|^q \Big),
$$
for  any $R>0$ such that $|x_1|\vee|x_2|\vee|y_1|\vee|y_2|\leq R,$ where the positive parameter $q\in(0,\frac{1}{2})$ the quantity $C_R$ depends on $R$ and $x\vee y$ denotes the maximum of $x, y.$(By the fact that we want the problem (\ref{ExSD-eq:general sde}) to be well-posed and by the conditions on $f_1,f_2$ and $g$ we get that $f_1,f_2,g$ are bounded on bounded intervals.)
\eass

Let the equidistant partition $0=t_0<t_1<...<t_N=T$ and $\D=T/N.$
We propose the following semi-discrete numerical scheme
\beqq\label{ExSD-eq:SD scheme} 
y_t=y_n + \int_{t_n}^{t} f_2(t_n, s, y_{t_n}, y_s)ds
+ \int_{t_n}^{t} g(t_n, s, y_{t_n}, y_s)dW_s,\quad t\in(t_n, t_{n+1}], 
\eeqq 
where we assume that for every $n\leq N-1,$ (\ref{ExSD-eq:SD scheme}) has a unique strong solution and $y_0=x_0$ a.s
$$y_n=y_{t_n} + f_1(t_n,y_{t_n})\cdot\D.$$
 In order to compare with the exact solution $x_t,$ which is a continuous time process, we consider the following interpolation process of the semi-discrete approximation, in a compact form, 
\beqq\label{ExSD-eq:SD scheme compact} 
y_t=y_0 + f_1(t_0,y_{t_0})\cdot\D + \int_{0}^{t}f_2(\hat{s}, s,y_{\hat{s}},y_s)ds +
\int_{0}^{t}g(\hat{s},s,y_{\hat{s}},y_s) dW_s, 
\eeqq 
where $\hat{s}=t_{n}$ when $s\in[t_n,t_{n+1}).$ Process (\ref{ExSD-eq:SD scheme compact}) has jumps at nodes $t_n.$ The first and third variable in $f_2,g$ denote the discretized part of the original SDE. We observe from (\ref{ExSD-eq:SD scheme compact}) that in order to solve for $y_t$, we have to solve an SDE and not an algebraic equation, thus in this context, we cannot reproduce implicit schemes, but we can reproduce the Euler scheme if we choose $f_1=0, f_2(s,r,x,y)=a(s,x)$ and $g(s,r,x,y)=b(s,x).$ 

The numerical scheme (\ref{ExSD-eq:SD scheme compact}) converges to the true solution $x_t$ of SDE (\ref{ExSD-eq:general sde}) and this is stated in the following, which  is  our main result.
\bth[Strong convergence]\label{ExSD-thm:strong_conv}
Suppose Assumption \ref{ExSD-ass:A} holds and (\ref{ExSD-eq:SD scheme}) has  a unique strong solution for every $n\leq N-1,$ where $x_0\in \bbl^p(\Omega,\bbR).$ Let also
$$
\bfE(\sup_{0\leq t\leq T}|x_t|^p) \vee \bfE(\sup_{0\leq t\leq T}|y_t|^p)<A,
$$
for some $p>2$ and $A>0.$ Then the semi-discrete numerical scheme (\ref{ExSD-eq:SD scheme compact}) converges to the true solution of (\ref{ExSD-eq:general sde}) in the $\bbl^2$-sense, that is
\beqq \label{ExSD-eq:strong_conv}
\lim_{\D\rightarrow0}\bfE\sup_{0\leq t\leq T}|y_t-x_t|^2=0.
\eeqq
\ethe

\section{Numerics}\label{NSF-sec:Experiments}
\setcounter{equation}{0}

Here, we make numerical tests to study the strong convergence of the proposed semi-discrete methods (\ref{NSF-eq:SD process}) and (\ref{NSF-eq:SD processNWalt}) for the Wright-Fisher model described by the It\^o SDE (\ref{NSF-eq:WFmodel}) with $k_1=A, k_2=A+B, k_3=C,$
where the parameters $A,B$ and $C$  are positive and $C=\sqrt{2k_2/(N_r-1)}.$ We take as initial condition the steady state of the deterministic part, that is $x_0=A/(A+B)$. This setting has been used for the approximation of ion channels within cardiac and neuronal cells, see \cite[Sec. 2.1]{dangerfield:2010}; the ion channel  occupies one of two positions (open and closed states) with transition rates $A$ and $B$ respectively and $N_r$ is the total number of ion channels within a cell,  see \cite[(2.3)]{dangerfield:2012}. We consider two set of parameters   
\begin{itemize}
\item SET I:   $(A, B, N_r)=(1, 2, 100),$
\item SET II:  $(A, B, N_r)=(7.0064, 0.0204, 100),$
\end{itemize}
as in \cite[Sec. 6 and 7.1]{dangerfield:2012}, where the Balance Implicit Split Step (BISS) method is suggested \cite[(4.8)]{dangerfield:2012} 
\beqq\label{NSF-eq:BISSscheme}
y_{n+1}^{BISS}=y_n + (A-(A+B)y_n)\D + \frac{C\sqrt{y_n(1-y_n)}\D W_n }{1+ d^1(y_n)|\D W_n |}(1-(A + B)\D), 
\eeqq
where $\D$ is the step-size of the equidistant discretization of the interval $[0,1]$, the control function $d^1$ is given by
$$
d^1(y) = \begin{cases}
    C\sqrt{(1-\vep)/\vep} \quad \,\text{if}\,\, y< \vep,\\        
    C\sqrt{(1-y)/y} \quad \text{if}\,\, \vep\leq y<1/2,\\
    C\sqrt{y/(1-y)} \quad \text{if}\,\, 1/2\leq y\leq 1-\vep,\\  
    C\sqrt{(1-\vep)/\vep} \quad \,\text{if}\,\, y>1 - \vep,
      \end{cases}
$$
and 
$$
\vep= \min \{ A\D, B\D, 1-A\D, 1-B\D\}.
$$
The hybrid (HYB)  scheme as proposed in \cite[(11)]{dangerfield:2010} is the result of a splitting method and reads
\beqq\label{NSF-eq:Hybscheme}
y_{n+1}^{HYB}=\frac{\alpha}{\beta}(e^{\beta\D}-1) + e^{\beta\D}\sin^2 \left(\frac{C}{2}\D W_n + \arcsin(\sqrt{y_n})\right).
\eeqq
It works only for the parameter SET I since we have to assume that 
$$\frac{a}{a+b}\in\left(\frac{1}{2(N_r-1)},1-\frac{1}{2(N_r-1)}\right).$$ 
Finally, the proposed semi-discrete (SD) scheme  reads
\beqq\label{NSF-eq:SDscheme1}
y_{n+1}^{SD}=\sin^2 \left(\frac{C}{2}\D W_n + \arcsin(\sqrt{y_n(1 + \beta \D) + \alpha\D})\right),
\eeqq
for SET I and 
\beqq\label{NSF-eq:SDscheme2}
y_{n+1}^{SD}=\sin^2 \left(\frac{C}{2}\D W_n + \arcsin\left(\sqrt{\frac{y_n(1 + \beta \D)+\alpha\D}{1+(\alpha+\beta)\D}}\right)\right),
\eeqq
for parameter SET II. The parameters of SET I are chosen in a way that the probability of the Euler-Maruyama (EM) scheme
$$  
y_{n+1}^{EM}= (A - (A+B)y_n)\D + C \sqrt{y_n(1-y_n)}\D W_n,
$$
leaving the interval $[0,1]$ is very small, whereas in the case of SET II this probability is high.  The paths of the solutions of EM exiting the boundaries $0$ and $1$  are only a few in the first case and one may reject them. Nevertheless, since such an approach induces bias to the solution obtained by the EM method (much more evident in the second case) we choose only to compare our method with BISS and HYB. 

We estimate the endpoint $\bbl^2$-norm $\ep=\sqrt{\bfE|y^{(\D)}(T) - x_T|^2},$ of the difference between the numerical scheme evaluated at step size $\D$ and the exact solution of (\ref{NSF-eq:WFmodel}). To do so, we compute $M$ batches of $L$ simulation paths, where each batch is estimated by 
$\hat{\ep}_j=\frac{1}{L}\sum_{i=1}^L|y_{i,j}^{(\D)}(T) - y_{i,j}^{(ref)}(T)|^2$ 
and the Monte Carlo estimator of the error is
\beqq\label{SF-def:l2error}
\hat{\ep}=\sqrt{\frac{1}{ML}\sum_{j=1}^M\sum_{i=1}^L|y_{i,j}^{(\D)}(T) - y_{i,j}^{(ref)}(T)|^2}
\eeqq
and requires $M\cdot L$ Monte Carlo sample paths. The reference solution is calculated using the method at a very fine time grid, $\D=2^{-13}.$ We have shown in Theorem \ref{NSF-theorem:Pol rate} and Proposition \ref{NSF-prop:Pol rate} that the SD numerical schemes converge strongly to the exact solution, so we use the SD method as a reference solution,
and the HYB method when applicable. The BISS method  converges in the $\bbl^1$-norm to the true solution \cite[Th. 5.1]{dangerfield:2012}, so we choose not to consider it as a reference solution even though we conjecture that a similar technique may be used to show an $\bbl^2$-convergence result.

\begin{figure}[th]
  \caption{ \small Convergence of SD, HYB and BISS methods applied to (\ref{NSF-eq:WFmodel}) with parameter SET I with HYB as a reference solution.} \label{fig:WFerrors_setI_HYB}
  \centering
   \includegraphics[width=0.7\textwidth]{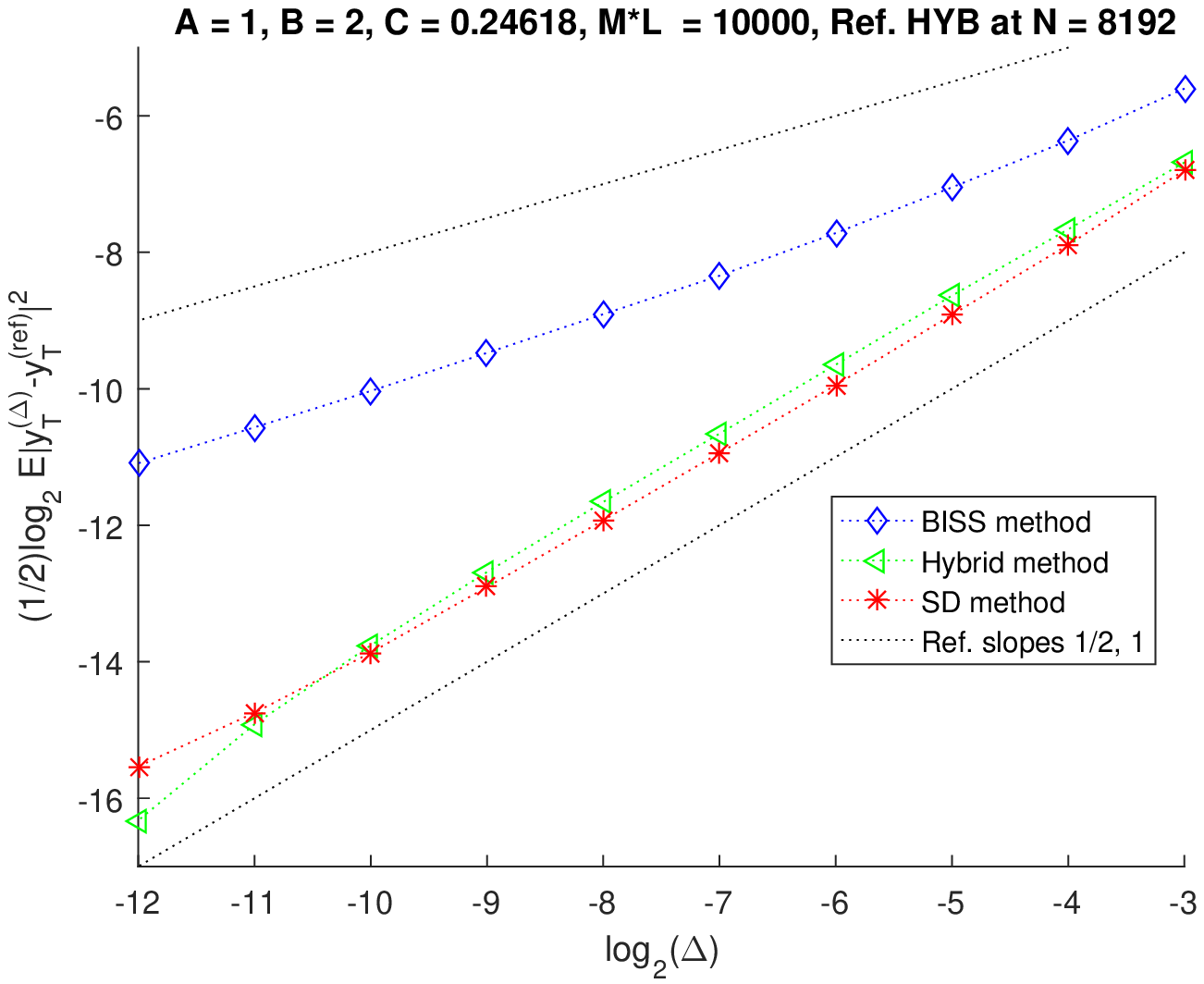} 
\end{figure}

\begin{figure}[th]
  \caption{ \small Convergence of SD, HYB and BISS methods applied to (\ref{NSF-eq:WFmodel}) with parameter SET I with SD as a reference solution.} \label{fig:WFerrors_setI_SD}
  \centering
   \includegraphics[width=0.7\textwidth]{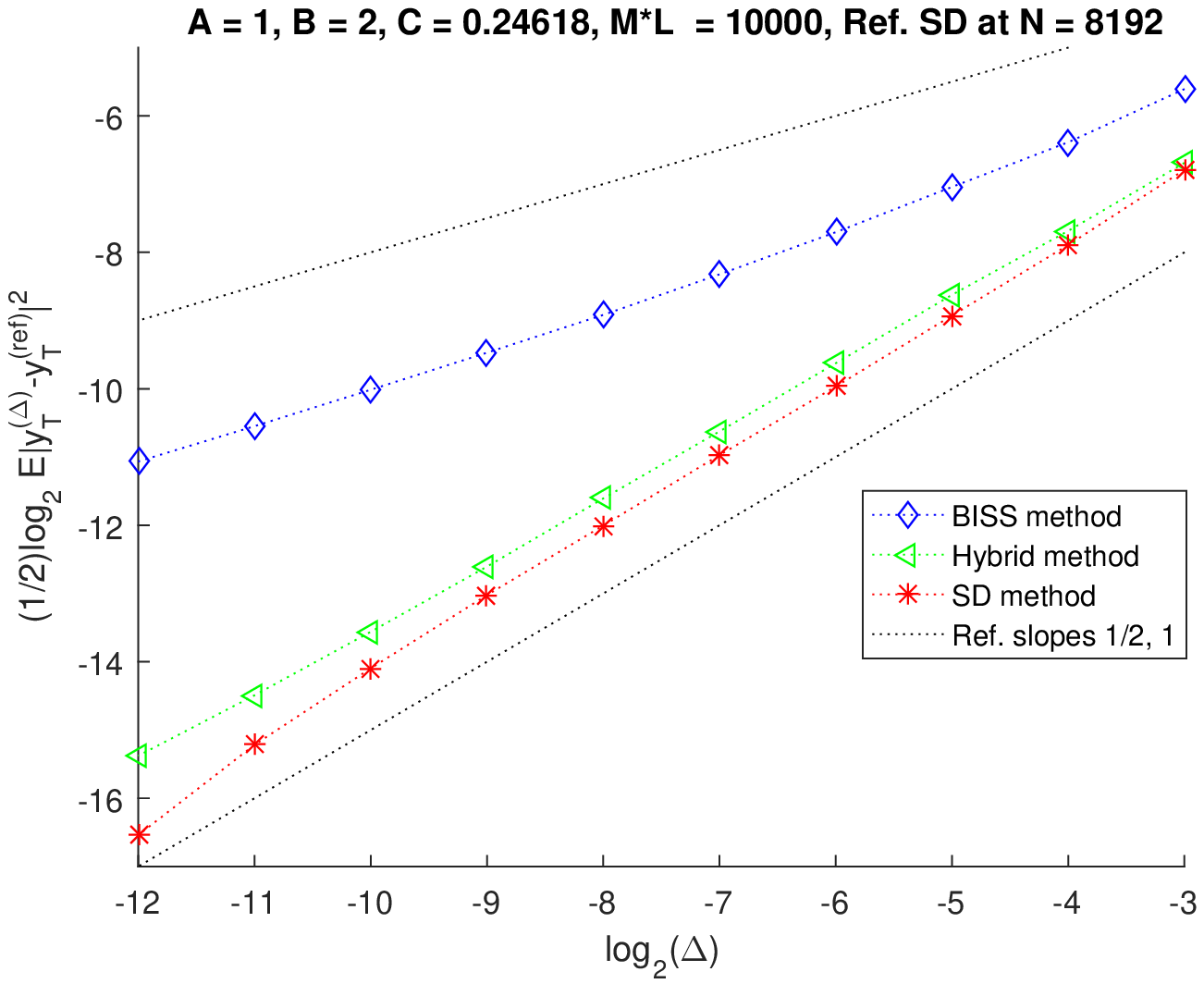} 
\end{figure}

\begin{figure}[th]
  \caption{ \small Convergence of SD and BISS methods applied to (\ref{NSF-eq:WFmodel}) with parameter SET II with SD as a reference solution.} \label{fig:WFerrors_setII_SD}
  \centering
   \includegraphics[width=0.7\textwidth]{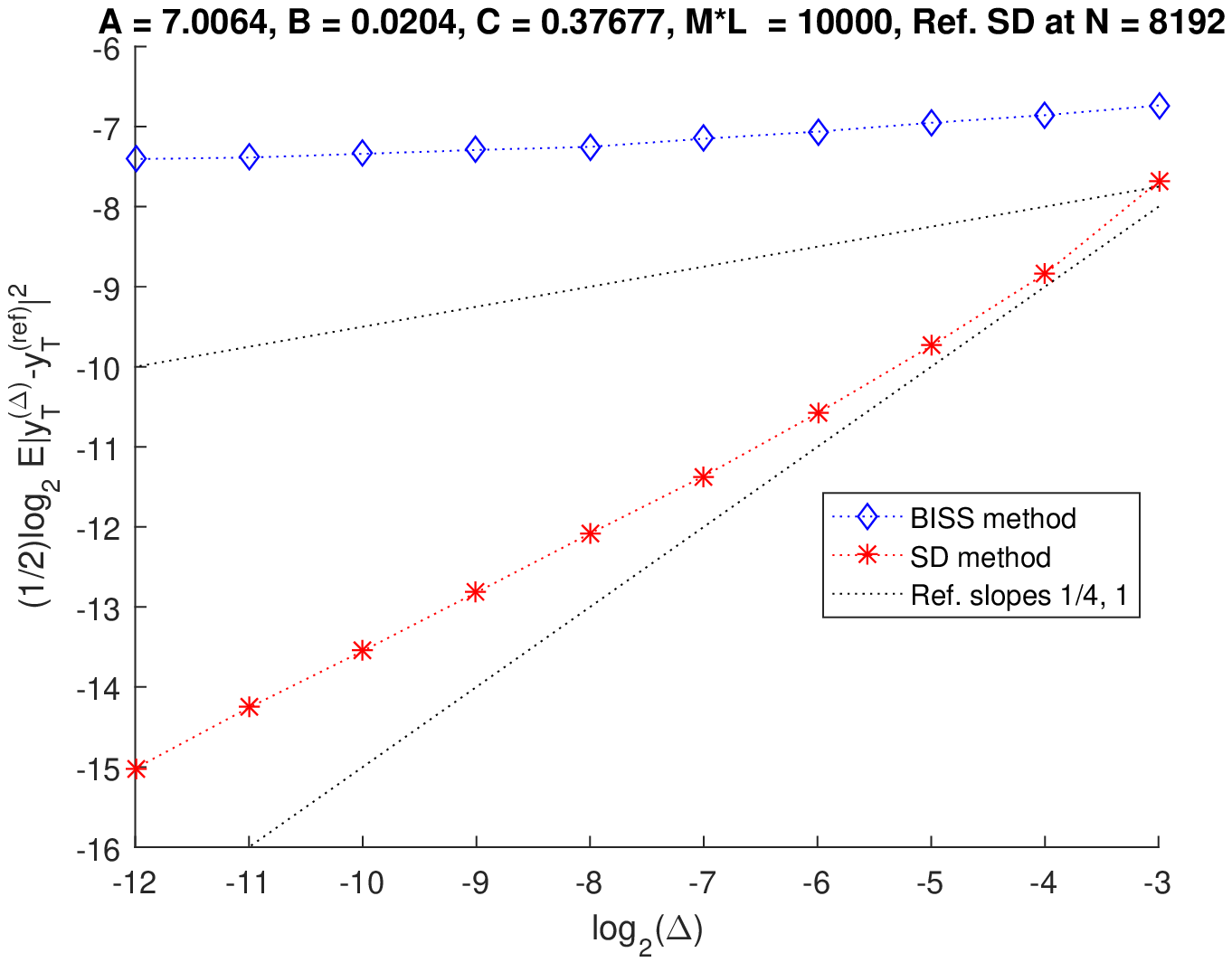} 
\end{figure}

We simulate $100\cdot 100=10^4$ paths, where the choice of the number of Monte Carlo paths  is adequately  large, so as not to significantly hinder the mean-square errors. We compute the approximation error (\ref{SF-def:l2error}) with $98\%$-confidence intervals.   
We  present the results for the parameter SET I and II in a $\log_2$-$\log_2$ scale in Figures \ref{fig:WFerrors_setI_HYB} - \ref{fig:WFerrors_setII_SD} and Tables \ref{tab:WFerrors_setI_with_HYB} - \ref{tab:WFerrors_setII_with_SD} respectively.

\begin{table}[htbp] 
\centering
        \begin{tabular}{c|r|r|r}
         Step  $\D$ & SD  & BISS  & HYB\\
 \hline $2^{-3}$  & $\textbf{0.008997}$       & $0.020613$     &     $0.009737$  \\
  \hline $2^{-4}$  & $\textbf{0.004216}$   & $0.012167$       &     $0.004937$      \\
 \hline $2^{-5}$  &  $\textbf{0.002072}$      &  $0.007558$    &    $0.002515$     \\
 \hline $2^{-6}$  &   $\textbf{0.001015}$      &  $0.004750$       &    $0.001248$      \\
 \hline $2^{-7}$  &   $\textbf{0.000509}$       &  $0.003080 $     &    $0.000619$     \\
 \hline $2^{-8}$  &  $\textbf{ 0.000257}$     &  $0.002084 $      &    $0.000309$    \\
 \hline $2^{-9}$  &  $\textbf{ 0.000130}$        &  $0.001407 $      &    $0.000152$  \\
 \hline $2^{-10}$  & $\textbf{0.000067}$    & $0.000955$        &    $0.000072$       \\
  \hline $2^{-11}$  & $0.000036$    & $0.000660$    &    $\textbf{0.000032}$   \\
   \hline $2^{-12}$  & $0.000021$     & $0.000459$     &    $\textbf{0.000012}$  
 \end{tabular}
    \caption{\small Error and step size of SD, BISS and Hybrid schemes for (\ref{NSF-eq:WFmodel}) with parameter SET I and HYB as a reference solution.}
    \label{tab:WFerrors_setI_with_HYB}  
\end{table}   

\begin{table}[htbp] 
\centering
        \begin{tabular}{c|r|r|r}
         Step  $\D$ & SD   & BISS  & HYB \\
 \hline $2^{-3}$  & $\textbf{0.009030}$       & $0.020507$         &   $0.009693$     \\
  \hline $2^{-4}$  & $\textbf{0.004210}$     & $0.011952$        & $0.004857$      \\
 \hline $2^{-5}$  &  $\textbf{0.002054}$       &  $0.007588$        &    $0.002536$        \\
 \hline $2^{-6}$  &   $\textbf{0.000999}$      &  $0.004791$        &    $0.001266$         \\
 \hline $2^{-7}$  &   $\textbf{0.000498}$      &  $0.003123 $       &    $0.000634$      \\
 \hline $2^{-8}$  &   $\textbf{0.000243}$       &  $0.002069 $       &     $0.000321$     \\
 \hline $2^{-9}$  &   $\textbf{0.000120}$       &  $0.001410 $       &     $0.000161$     \\
 \hline $2^{-10}$  & $\textbf{0.000057}$       & $0.000967$        &    $0.000083$      \\
  \hline $2^{-11}$  & $\textbf{0.000026}$      & $0.000668$       &   $0.000043$   \\
   \hline $2^{-12}$  & $\textbf{0.000011}$      & $0.000467$       &    $0.000023$      
 \end{tabular}
    \caption{\small Error and step size of SD, BISS and Hybrid schemes for (\ref{NSF-eq:WFmodel}) with parameter SET I and SD as a reference solution.}
    \label{tab:WFerrors_setI_with_SD}  
\end{table}

\begin{table}[htbp] 
\centering
        \begin{tabular}{c|r|r}
         Step  $\D$ & SD   & BISS \\
 \hline $2^{-3}$  & $\textbf{0.004860}$               &   $0.009378$     \\
  \hline $2^{-4}$  & $\textbf{0.002181}$             & $0.008626$      \\
 \hline $2^{-5}$  &  $\textbf{0.001169}$          &    $0.008074$        \\
 \hline $2^{-6}$  &   $\textbf{0.000653}$      &  $0.007477$          \\
 \hline $2^{-7}$  &   $\textbf{0.000376}$      &  $0.007035 $        \\
 \hline $2^{-8}$  &   $\textbf{0.000231}$       &  $0.006556 $         \\
 \hline $2^{-9}$  &   $\textbf{0.000138}$       &  $0.006381 $          \\
 \hline $2^{-10}$  & $\textbf{0.000083}$       & $0.006171$             \\
  \hline $2^{-11}$  & $\textbf{0.000051}$      & $0.005977$       \\
   \hline $2^{-12}$  & $\textbf{0.000030}$      & $0.005898$       
 \end{tabular}
    \caption{\small Error and step size of SD and BISS schemes for (\ref{NSF-eq:WFmodel}) with parameter SET II and SD as a reference solution.}
    \label{tab:WFerrors_setII_with_SD}
\end{table}

\begin{figure}[th]
  \caption{ \small Convergence of SD, EM and BISS methods applied to (\ref{NSF-eq:WFmultimodel})-(\ref{NSF-eq:WFmultimodel2}) with parameter SET III with EM as a reference solution.} \label{fig:WFerrors_setIII_EM}
  \centering
   \includegraphics[width=0.7\textwidth]{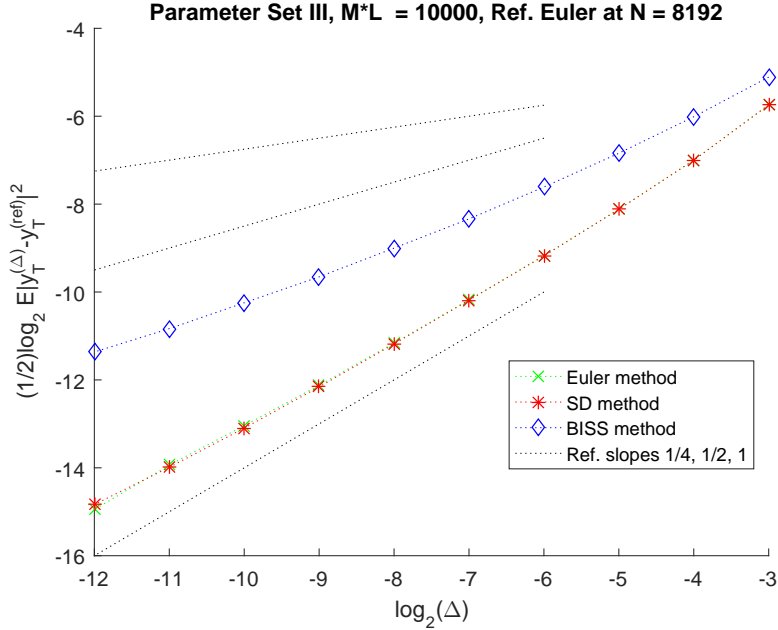} 
\end{figure}

Furthermore, we  study the strong convergence of the proposed semi-discrete methods (\ref{NSF-eq:SDmultiprocess1})-(\ref{NSF-eq:SDmultiprocess2}) for the Wright-Fisher model described by the  system of It\^o SDEs (\ref{NSF-eq:WFmultimodel})-(\ref{NSF-eq:WFmultimodel2}) with $k_1^{(1,1)}=A_3, k_1^{(1,2)}=A_2-A_3, k_1^{(2,1)}=A_1, k_1^{(2,2)}=B_1-A_1, k_2^{(1)}=B_3+B_1+A_3, k_2^{(2)}=A_2+B_2+A_1, k_3^{(1,1)}=-C_1, k_3^{(1,2)}=C_2, k_3^{(2,1)}=C_1, k_3^{(2,3)}=-C_3.$ This setting has been used for the approximation of the proportion of alleles or channels in a 3 state system and we consider the set of parameters   
\begin{itemize}
\item SET III:   $(A_1, A_2, A_3, B_1, B_2, B_3, C_1, C_2, C_3)=(1, 2, 3, 1.2, 2.3, 3.4, 0.1271, 0.1798, 0.1291),$
\end{itemize}
as in \cite[Sec. 6]{dangerfield:2012}, where the Balance Implicit Split Step (BISS) method is suggested; they use the splitting method, solving first the stochastic system with a balanced implicit scheme 
\beam\nonumber
_{BISS}Y_{n+1}^{(1)}&=&\frac{Y_n^{(1)} - C_1\sqrt{Y_n^{(1)}Y_n^{(2)}}\D W_n^{(1)} + C_2\sqrt{Y_n^{(1)}(1-Y_n^{(1)}-Y_n^{(2)})}\D W_n^{(2)}}{1+D_1(Y_n) + D_2(Y_n) + D_3(Y_n)}\\
\label{NSF-eq:BISSschememulti1}&& + \frac{Y_n^{(1)}(1+D_1(Y_n) + D_2(Y_n) + D_3(Y_n))}{1+D_1(Y_n) + D_2(Y_n) + D_3(Y_n)}, 
\eeam
\beam\nonumber
_{BISS}Y_{n+1}^{(2)}&=&\frac{Y_n^{(2)} + C_1\sqrt{Y_n^{(1)}Y_n^{(2)}}\D W_n^{(1)} - C_3\sqrt{Y_n^{(2)}(1-Y_n^{(1)}-Y_n^{(2)})}\D W_n^{(3)}}{1+D_1(Y_n) + D_2(Y_n) + D_3(Y_n)}\\
\label{NSF-eq:BISSschememulti2}&& + \frac{Y_n^{(2)}(1+D_1(Y_n) + D_2(Y_n) + D_3(Y_n))}{1+D_1(Y_n) + D_2(Y_n) + D_3(Y_n)}, 
\eeam
and then the deterministic part with the one-step EM scheme; here the control functions $D_i, i=1,\ldots, 3,$ are given by
$$
D_1(Y) = \begin{cases}
    C_1\left(\sqrt{\frac{Y^{(2)}}{Y^{(1)}}} +  \sqrt{\frac{\vep}{Y^{(1)}Y^{(2)}}}\right)|\D W_n^{(1)}|\quad \,\text{if}\,\, \vep<Y^{(1)}\leq Y^{(2)},\\        
    C_1\left(\sqrt{\frac{1-Y^{(1)}-Y^{(2)}}{Y^{(2)}}} +  \sqrt{\frac{\vep}{Y^{(1)}Y^{(2)}}}\right)|\D W_n^{(1)}|\quad \,\text{if}\,\, \vep<Y^{(2)}< Y^{(1)},
      \end{cases}
$$
$$
D_2(Y) = \begin{cases}
    C_2\left(\sqrt{\frac{1-Y^{(1)}-Y^{(2)}}{Y^{(1)}}} +  \sqrt{\frac{\vep}{Y^{(1)}(1-Y^{(1)}-Y^{(2)})}}\right)|\D W_n^{(2)}|\quad \,\text{if}\,\, 2Y^{(1)}+Y^{(2)}<1,\\        
    C_2\left(\sqrt{\frac{Y^{(1)}}{1-Y^{(1)}-Y^{(2)}}} +  \sqrt{\frac{\vep}{Y^{(1)}(1-Y^{(1)}-Y^{(2)})}}\right)|\D W_n^{(2)}|\quad \,\text{if}\,\, 2Y^{(1)}+ Y^{(2)}\geq1,
      \end{cases}
$$
$$
D_3(Y) = \begin{cases}
    C_3\left(\sqrt{\frac{1-Y^{(1)}-Y^{(2)}}{Y^{(2)}}} +  \sqrt{\frac{\vep}{Y^{(2)}(1-Y^{(1)}-Y^{(2)})}}\right)|\D W_n^{(3)}|\quad \,\text{if}\,\, 2Y^{(2)}+Y^{(1)}<1,\\        
    C_3\left(\sqrt{\frac{Y^{(2)}}{1-Y^{(1)}-Y^{(2)}}} +  \sqrt{\frac{\vep}{Y^{(2)}(1-Y^{(1)}-Y^{(2)})}}\right)|\D W_n^{(3)}|\quad \,\text{if}\,\, 2Y^{(2)}+ Y^{(1)}\geq1;
      \end{cases}
$$ 
In case $Y^{(i)}\leq\vep$ then set $Y^{(i)}=\vep, i=1,2$ and if $Y^{(1)}+Y^{(2)}\geq1-\vep$ then set $Y^{(1)}+Y^{(2)}=1-\vep.$ The proposed semi-discrete scheme reads, see (\ref{NSF-eq:SDmultisolution1})-(\ref{NSF-eq:SDmultisolution2})
\beqq\label{NSF-eq:SDmultischeme1}
_{SD}Y_{n+1}^{(1)}=\sin^2 \left(-\frac{C_1}{2}\sqrt{\frac{Y_{n}^{(2)}}{1-Y_{n}^{(1)}}}\D W_n^{(1)} + \frac{C_2}{2}\sqrt{\frac{1-Y_{n}^{(1)}-Y_{n}^{(2)}}{1-Y_{n}^{(1)}}}\D W_n^{(2)}  + \arcsin(\sqrt{y_n^{(1)}})\right),
\eeqq
\beqq\label{NSF-eq:SDmultischeme2}
_{SD}Y_{n+1}^{(2)} = 
\sin \left(\frac{C_1}{2}\sqrt{\frac{Y_{n}^{(1)}}{1-Y_{n}^{(2)}}}\D W_n^{(1)} - \frac{C_3}{2}\sqrt{\frac{1-Y_{n}^{(1)}-Y_{n}^{(2)}}{1-Y_{n}^{(2)}}}\D W_n^{(3)} + \arcsin(\sqrt{y_n^{(2)}})\right),
\eeqq
where $\D$ is the step-size of the equidistant discretization of the interval $[0,1].$ Note that for the parameter SET III $0<y_n^{(i)}<1, i=1,2,$ for any $\D<\min\{1/k_2^{(1)},1/k_2^{(2)}\}= 1/(B_3+B_1+A_3).$
The parameters of SET III are chosen in a way that the probability of the Euler-Maruyama (EM) scheme
\beam 
\nonumber 
_{EM}Y_{n+1}^{(1)} &=& Y_n^{(1)} + \left(A_3+(A_2-A_3)Y_n^{(2)} -(B_3+B_1+A_3)Y_n^{(1)}\right)\D\\
\label{NSF-eq:EMmultischeme1}&&- C_1\sqrt{Y_n^{(1)}Y_n^{(2)}}\D W_n^{(1)} + C_2\sqrt{Y_n^{(1)}(1-Y_n^{(1)}-Y_n^{(2)})}\D W_n^{(2)},\\
\nonumber 
_{EM}Y_{n+1}^{(2)} &=& Y_n^{(2)} + \left(A_1+(B_1-A_1)Y_n^{(1)} -(A_2+B_2+A_1)Y_n^{(2)}\right)\D\\
\label{NSF-eq:EMmultischeme1}&&+ C_1\sqrt{Y_n^{(1)}Y_n^{(2)}}\D W_n^{(1)} - C_3\sqrt{Y_n^{(2)}(1-Y_n^{(1)}-Y_n^{(2)})}\D W_n^{(3)},
\eeam
leaving the region $[0,1]^3$ is very small.  The paths of the solutions of EM exiting the boundaries $0$ and $1$  are extremely few when we take as initial condition the steady state of the deterministic part of the system (which for the parameter SET III is away from the boundaries) and consider a very fine time discretization ($\D=2^{-13}$). Thus, we take as the exact solution the EM method and reject the paths, if any, outside the region $[0,1]^3.$ We compare our method with BISS and present the results in  Figure \ref{fig:WFerrors_setIII_EM}.

We make the  following comments.
\begin{itemize}
    \item The performance of the HYB and SD schemes for parameter Set I  is quite similar with SD producing smaller errors. They both perform better than BISS. Nevertheless, HYB is not boundary preserving for parameter  SET II. The performance of EM and SD schemes is almost identical for parameter Set III. They both perform better than BISS. Nevertheless, EM is not boundary preserving (we just used it in this experiment for comparative reasons as in \cite[Sec. 6]{dangerfield:2012}.)
    \item The numerical results suggest that the SD schemes converge in the mean-square sense with order $1$ for parameter SET I and SET III and at least $1/2$ for parameter SET II. 
   \item The proposed  SD schemes perform better  w.r.t. the computational time required to achieve a desired level of accuracy, since there is no need to calculate a control function.
    \item For the implementation of the SD method, we have to assume that $\D<-1/(\alpha +\beta)=0.1423,$ for the parameters of SET I, $\D<-1/\beta=0.1437,$ for the parameters of SET II,  and $\D<1/(B_1+B_3+A_3)=1/7.6,$ for the parameters of SET III, so the step-size $\D=1/8$ is sufficient. There is not a step-size restriction in the BISS method; nevertheless we propose the SD method since it is fast and more accurate.  
\end{itemize}


\section{Proof of Theorem \ref{ExSD-thm:strong_conv}.}\label{ExSD-sec:Proof}
\setcounter{equation}{0}

We split the proof is three steps. First, we prove a general estimate of the error of the SD method for any $p>0.$ Then, we show the $\bbl^1$-convergence of the semi-discrete method and finally the  $\bbl^2$-convergence (\ref{ExSD-eq:strong_conv}).
We denote the indicator function of a set $A$ by $\bbi_{A}.$ The quantity $C_R$ may vary from line to line and it may depend apart from $R$ on other quantities, like time $T$ for example, which are all constant, in the sense that we don't let them grow to infinity.

\ble[Error bound for the semi-discrete scheme] \label{ExSD-lem:Error_SDbound} 
Let the assumption of Theorem \ref{ExSD-thm:strong_conv} hold. Let $R>0,$ and set the stopping time $\theta_R=\inf\{t\in [0,T]: |y_t|>R \, \hbox { or } \,  |y_{\hat{t}}|>R\}.$ Then the following estimate holds 
$$ 
\bfE|y_{s\wedge\theta_R}-y_{\wh{s\wedge\theta_R}}|^p \leq C_{R}\D^{p/2}, $$ 
for any $p>0,$ where  $C_R$ does not depend on $\D,$ implying
$\sup_{s\in[t_{n_s},t_{n_s+1}]}\bfE|y_{s\wedge\theta_R}-y_{\wh{s\wedge\theta_R}}|^p=O(\D^{p/2}),$ as $\D\downarrow0.$ 
\ele

\bpf[Proof of Lemma \ref{ExSD-lem:Error_SDbound}]
We fix a $p\geq2.$ Let $n_s$ integer such that $s\in[t_{n_s},t_{n_s+1}).$ It holds that
\beao
&&|y_{s\wedge\theta_R}-y_{\wh{s\wedge\theta_R}}|^p=\left| \int_{t_{\wh{n_s\wedge\theta_R}}}^{s\wedge\theta_R}f_2(\hat{u},u,y_{\hat{u}},y_u)du + \int_{t_{\wh{n_s\wedge\theta_R}}}^{s\wedge\theta_R}g(\hat{u},u,y_{\hat{u}},y_u) dW_u\right|^p\\
&\leq&2^{p-1}\left|\int_{t_{\wh{n_s\wedge\theta_R}}}^{s\wedge\theta_R}f_2(\hat{u},u,y_{\hat{u}},y_u)du\right|^p + 2^{p-1}\left|\int_{t_{\wh{n_s\wedge\theta_R}}}^{s\wedge\theta_R}g(\hat{u},u,y_{\hat{u}},y_u) dW_u\right|^p\\
&\leq&2^{p-1}|s\wedge\theta_R-t_{\wh{n_s\wedge\theta_R}}|^{p-1}\int_{t_{\wh{n_s\wedge\theta_R}}}^{s\wedge\theta_R}|f_2(\hat{u},u,y_{\hat{u}},y_u)|^pdu + 2^{p-1}\left|\int_{t_{\wh{n_s\wedge\theta_R}}}^{s\wedge\theta_R}g(\hat{u},u,y_{\hat{u}},y_u) dW_u\right|^p\\
&\leq&C_R\D^{p} + 2^{p-1}\left|\int_{t_{\wh{n_s\wedge\theta_R}}}^{s\wedge\theta_R}g(\hat{u},u,y_{\hat{u}},y_u)dW_u\right|^p, \eeao 
where we have used Cauchy-Schwarz inequality and Assumption \ref{ExSD-ass:A} for the function $f_2$.
 Taking expectations  in the above inequality gives
\beao
\bfE|y_{s\wedge\theta_R}-y_{\wh{s\wedge\theta_R}}|^p &\leq& C_R\D^{p}  + 2^{p-1}\bfE\left|\int_{t_{\wh{n_s\wedge\theta_R}}}^{ t_{n_s+1}\wedge\theta_R }g(\hat{u},u,y_{\hat{u}},y_u) dW_u\right|^p\\
&\leq& C_R\D^{p}  +  2^{p-1}\un{\left(\frac{p^{p+1}}{2(p-1)^{p-1}}\right)^{p/2}}_{C_p}\bfE\left|\int_{t_{\wh{n_s\wedge\theta_R}}}^{ t_{n_s+1}\wedge\theta_R }|g(\hat{u},u,y_{\hat{u}},,y_u) |^2du\right|^{p/2}\\
&\leq& C_R\D^{p}  +  2^{p-1}C_p\D^{\frac{p-2}{2}}\bfE\int_{t_{\wh{n_s\wedge\theta_R}}}^{ t_{n_s+1}\wedge\theta_R }|g(\hat{u},u,y_{\hat{u}},y_u) |^p du\\
&\leq& C_R\D^{p}  +  C_R\D^{p/2},
\eeao
where in the third step we have used the Burkholder-Davis-Gundy (BDG)  inequality \cite[Th. 1.7.3]{mao:2007}, \cite[Th. 3.3.28]{karatzas_shreve:1988} on the diffusion term,  in the last step Assumption \ref{ExSD-ass:A} for the function $g.$  Thus,
$$
\lim_{\D\downarrow0}\frac{\sup_{s\in[t_{n_s},t_{n_s+1}]}\bfE|y_{s\wedge\theta_R}-y_{\wh{s\wedge\theta_R}}|^p}{\D^{p/2}}\leq C_R,
$$
which justifies the $O(\D^{p/2})$ notation.   
Now for $0<p<2$ we have that
$$
\bfE|y_{s\wedge\theta_R}-y_{\wh{s\wedge\theta_R}}|^p\leq\left(\bfE|y_{s\wedge\theta_R}-y_{\wh{s\wedge\theta_R}}|^2\right)^{p/2}\leq C_R\D^{p/2},
$$
where we have used Jensen inequality for the concave function $\phi(x)=x^{p/2}.$ 
\epf

In the next result we estimate the $\bbl^1$-error using the Yamada-Watanabe approach. We denote the difference $\bbE_{t}:=y_{t}-x_{t},$

\bpr\label{ExSD-prop:L1_conv}[Convergence of the semi-discrete scheme in $\bbl^1$]
Let the assumptions of Theorem \ref{ExSD-thm:strong_conv} hold. Let $R>0,$ and set the stopping time $\theta_R=\inf\{t\in [0,T]: |y_t|>R \, \hbox { or } \,  |x_{t}|>R\}.$ Then we have
\beqq \label{ExSD-eq:L1_conv}
\sup_{0\leq t\leq T}\bfE|\bbE_{t\wedge\theta_R}|\leq \left( C_R\sqrt{\D} + \frac{C_R}{me_m}\D^{q} +\frac{C_R}{m}e_m +  e_{m-1}\right)e^{a_{R,m,q}t},
\eeqq
for any $m>1,$ where
$$
e_m=e^{-m(m+1)/2}, \quad a_{R,m,q}:=C_R + \frac{C_R}{m} +  \frac{C_R}{m(e_m)^{\frac{2-2q}{q}}}
$$
and $C_R$ does not depend on $\D.$ It holds that $\lim_{m\uparrow\infty}e_m=0.$
\epr
\bpf[Proof of Proposition \ref{ExSD-prop:L1_conv}]
Let the non-increasing sequence $\{e_m\}_{m\in\bbN}$ with $e_m=e^{-m(m+1)/2}$ and $e_0=1.$ We introduce the following sequence of smooth approximations of $|x|,$ (method of Yamada and Watanabe, \cite{yamada_watanabe:1971})
$$
\phi_m(x)=\int_0^{|x|}dy\int_0^{y}\psi_m(u)du,
$$
where the existence of the continuous function $\psi_m(u)$ with $0\leq \psi_m(u) \leq 2/(mu)$ and support in $(e_m,e_{m-1})$ is justified by $\int_{e_m}^{e_{m-1}}(du/u)=m.$ The following relations hold for $\phi_m\in\bbc^2(\bbR,\bbR)$ with $\phi_m(0)=0,$
 $$
 |x| - e_{m-1}\leq\phi_m(x)\leq |x|, \quad |\phi_{m}^{\prime}(x)|\leq1, \quad x\in\bbR, $$
 $$
 |\phi_{m}^{\prime \prime }(x)|\leq\frac{2}{m|x|}, \,\hbox{ when }  \,e_m<|x|<e_{m-1} \,\hbox{ and }  \,  |\phi_{m}^{\prime \prime }(x)|=0 \,\hbox{ otherwise. }
 $$
We have that
\beqq\label{ExSD-eq:YW}
\bfE|\bbE_{t\wedge\theta_R}| \leq e_{m-1} + \bfE\phi_m(\bbE_{t\wedge\theta_R}).
\eeqq
Applying It\^o's formula to the sequence $\{\phi_m\}_{m\in\bbN},$ we get
\beao
&&\phi_m(\bbE_{t\wedge\theta_R})=\int_{0}^{t\wedge\theta_R} \phi_m^{\prime}(\bbE_s) (f_2(\hat{s},s,y_{\hat{s}},y_s)-f_2(s,s,x_s,x_s)-f_1(s,x))ds +M_t\\
&&+ \frac{1}{2}\int_{0}^{t\wedge\theta_R} \phi_m^{\prime\prime}(\bbE_s) (g(\hat{s},s,y_{\hat{s}},y_s)-g(s,s,x_s,x_s))^2 ds \\
\nonumber&\leq& \int_{0}^{t\wedge\theta_R} C_R\left( |y_{\hat{s}}-x_s | + |\bbE_s| + |\hat{s}-s|  \right) ds + M_t\\
&& + C_R\int_{0}^{t\wedge\theta_R} \frac{1}{m|\bbE_s|} \Big( |y_s - y_{\hat{s}}|^2 + |\bbE_s|^2  + |\hat{s}-s|^2 + \left(|y_{\hat{s}}-y_s|^{2q} + |\bbE_s|^{2q}\right)\Big)ds\\
&\leq& C_R\int_{0}^{t\wedge\theta_R}|y_s-y_{\hat{s}}|ds + C_R \int_{0}^{t\wedge\theta_R} |\bbE_s|ds +  C_R \int_{0}^{t\wedge\theta_R} |\hat{s}-s|ds + M_t\\
&&+  \frac{C_R}{m}\int_{0}^{t\wedge\theta_R}\frac{|y_s - y_{\hat{s}}|^2 + |y_{\hat{s}}-y_s|^{2q} + |\bbE_s|^2 + |\bbE_s|^{2q} +|\hat{s}-s|^2}{|\bbE_s|}ds\\
&\leq& C_R\int_{0}^{t\wedge\theta_R}|y_s-y_{\hat{s}}|ds +\frac{C_R}{me_m}\int_{0}^{t\wedge\theta_R}\left(|y_s-y_{\hat{s}}|^2 +|y_s-y_{\hat{s}}|^{2q}\right)ds + (C_R + \frac{C_R}{m}) \int_{0}^{t\wedge\theta_R} |\bbE_s|ds\\
&& + \frac{C_R}{m}\int_{0}^{t\wedge\theta_R} |\bbE_s|^{2q-1}ds +  \frac{C_R}{me_m}\sum_{k=0}^{[t/\D-1]}\int_{t_k}^{t_{k+1}\wedge\theta_R}|t_k-s|^2ds + C_R\sum_{k=0}^{[t/\D-1]}\int_{t_k}^{t_{k+1}\wedge\theta_R}|t_k-s|ds + M_t,
\eeao
where in the second step we have used Assumption \ref{ExSD-ass:A} for the functions $f_1,f_2, g$ the subadditivity property of $h(x)=x^{2q},$ and the properties of $\phi_m$ and
$$
M_t:= \int_{0}^{t\wedge\theta_R} \phi_m^{\prime}(\bbE_u)(g(\hat{u},u,y_{\hat{u}},y_u)-g(u,u,x_u,x_u)) dW_u.
$$
Using the estimate 
\beao
\frac{C_R}{m}\int_{0}^{t\wedge\theta_R} |\bbE_s|^{2q-1}ds&\leq&\frac{C_R}{m|\bbE_s|}\int_{0}^{t\wedge\theta_R} \left(q|\bbE_s|^{2}(e_m)^{\frac{2q-2}{q}} + (1-q)(e_m)^2 \right)ds\\
&\leq&\frac{C_R}{m(e_m)^{\frac{2-2q}{q}}}\int_{0}^{t\wedge\theta_R}|\bbE_s|ds + \frac{C_R}{m}e_m,
\eeao
we get
\beao
&&\phi_m(\bbE_{t\wedge\theta_R})\leq C_R\int_{0}^{t\wedge\theta_R}|y_s-y_{\hat{s}}|ds +\frac{C_R}{me_m}\int_{0}^{t\wedge\theta_R}\left(|y_s-y_{\hat{s}}|^2 +|y_s-y_{\hat{s}}|^{2q}\right)ds\\
&& + \left(C_R + \frac{C_R}{m} +  \frac{C_R}{m(e_m)^{\frac{2-2q}{q}}}\right) \int_{0}^{t\wedge\theta_R} |\bbE_s|ds + \frac{C_R}{m}e_m +  \frac{C_R}{me_m}\D^2 + C_R\D + M_t.
\eeao
Taking expectations in the above inequality yields
$$\bfE\phi_m(\bbE_{t\wedge\theta_R})\leq C_R\D + C_R\sqrt{\D} + \frac{C_R}{me_m}(\D^2 + \D + \D^{q}) +\frac{C_R}{m}e_m
+ \left(C_R + \frac{C_R}{m} +  \frac{C_R}{m(e_m)^{\frac{2-2q}{q}}}\right) \int_{0}^{t\wedge\theta_R} \bfE|\bbE_s|ds,
$$
where we have used  Lemma \ref{ExSD-lem:Error_SDbound} and the fact that $\bfE M_t=0$.\footnote{The function $h(u)=\phi_m^{\prime}(\bbE_u)(g(\hat{u},u,y_{\hat{u}},y_u)-g(u,u,x_u,x_u))$ belongs to the space $\bbm^2([0,t\wedge\theta_R];\bbR)$ of real-valued measurable $\bbf_t$-adapted processes such that $\bfE\int_0^{t\wedge\theta_R}|h(u)|^2du<\infty$ thus (\cite[Th. 1.5.8]{mao:2007}) implies $\bfE M_t=0$.} Thus (\ref{ExSD-eq:YW}) becomes
\beao
\bfE|\bbE_{t\wedge\theta_R}|&\leq& C_R\sqrt{\D} + \frac{C_R}{me_m}\D^{q} +\frac{C_R}{m}e_m +  e_{m-1} + \left(C_R + \frac{C_R}{m} +  \frac{C_R}{m(e_m)^{\frac{2-2q}{q}}}\right) \int_{0}^{t\wedge\theta_R} \bfE|\bbE_s|ds\\
&\leq&\left( C_R\sqrt{\D} + \frac{C_R}{me_m}\D^{q} +\frac{C_R}{m}e_m +  e_{m-1}\right)e^{a_{R,m,q}t},
\eeao
where in the last step we have used the Gronwall inequality (\cite[(7)]{gronwall:1919}) and $a_{R,m,q}=C_R + \frac{C_R}{m} +  \frac{C_R}{m(e_m)^{\frac{2-2q}{q}}}.$ Taking the supremum over all $0\leq t \leq T$ gives (\ref{ExSD-eq:L1_conv}).
\epf

\subsection*{Convergence of the semi-discrete scheme in $\bbl^2$.}

Let the events $\W$ be defined by $\W_R:=\{\w \in\W: \sup_{0\leq t\leq T}|x_t|\leq R, \sup_{0\leq t\leq T}|y_t|\leq R\}$ and the stopping time $\theta_R=\{\inf t\in[0,T]: |y_t|>R \, \hbox{ or } \,|x_t|>R\}$ for some $R>0$ big enough. We have that
\beam
\nonumber
\bfE\sup_{0\leq t\leq T}|\bbE_t|^2 &=& \bfE\sup_{0\leq t\leq T}|\bbE_t|^2\bbi_{\W_R} +
\bfE\sup_{0\leq t\leq T}|\bbE_t|^2\bbi_{(\W_R)^c}\\
\nonumber&\leq& \bfE\sup_{0\leq t\leq T }|\bbE_{t\wedge_R}|^2 +
\left(\bfE\sup_{0\leq t\leq T}|\bbE_t|^p\right)^{2/p}\left(\bfE(\bbi_{(\W_R)^c})^{2p/(p-2)}\right)^{(p-2)/p}\\
\nonumber&\leq& \bfE\sup_{0\leq t\leq T }|\bbE_{t\wedge\theta_R}|^2 + 
\left(\bfE\sup_{0\leq t\leq T}|\bbE_t|^p\right)^{2/p}\left(\bfP(\W_R)^c\right)^{(p-2)/p}\\
\nonumber&\leq& \bfE\sup_{0\leq t\leq T }|\bbE_{t\wedge\theta_R}|^2 + 
\left(2^{p-1}\bfE\sup_{0\leq t\leq T}(|y_t|^p + |x_t|^p)\right)^{2/p}\left(\bfP(\W_R)^c\right)^{(p-2)/p}\\
\label{ExSD-eq:L2_conv}&\leq& \bfE\sup_{0\leq t\leq T }|\bbE_{t\wedge\theta_R}|^2 + 4\cdot A^{2/p}\left(\bfP(\W_R)^c\right)^{(p-2)/p},
\eeam
where $p>2$ is such that the moments of $|x_t|^p$ and $|y_t|^p$ are bounded by the constant $A.$ We want to estimate each term of the right hand side of (\ref{ExSD-eq:L2_conv}). It holds that
\beao
\bfP(\W_R^c) &\leq& \bfP(\sup_{0\leq t\leq T}|y_t|>R) +  \bfP(\sup_{0\leq t\leq T}|x_t|>R)\\
&\leq& (\bfE\sup_{0\leq t\leq T}|y_t|^k)R^{-k} +  (\bfE\sup_{0\leq t\leq T}|x_t|^k)R^{-k},
\eeao
for any $k\geq1$ where the first step comes from the subadditivity of the measure $\bfP$ and the second step from Markov inequality. Thus for $k=p$ we get 
$$
\bfP(\W_R^c)\leq 2AR^{-p}.
$$
We estimate the difference $|\bbE_{t\wedge\theta_R}|^2=|y_{t\wedge\theta_R}-x_{t\wedge\theta_R}|^2.$ It\^o's formula implies that
\beao
&&|\bbE_{t\wedge\theta_R}|^2=\int_{0}^{t\wedge\theta_R}2\left(f_2(\hat{s},s,y_{\hat{s}},y_s)-f_2(s,s,x_s,x_s)-f_1(s,x_s)\right)|\bbE_s| + \left( g(\hat{s},s,y_{\hat{s}},y_s)-g(s,s,x_s,x_s) \right)^2 ds\\
&&+ |f_1(t_0,y_{t_0})\D|^2 + \int_{0}^{t\wedge\theta_R}2|\bbE_s|\left( g(\hat{s},s,y_{\hat{s}},y_s)-g(s,s,x_s,x_s)\right)dW_s\\
&\leq&\int_{0}^{t\wedge\theta_R}|f_2(\hat{s},s,y_{\hat{s}},y_s)-f_2(s,s,x_s,x_s)|^2ds + \int_{0}^{t\wedge\theta_R} |\bbE_s|^2 ds + 
2\int_{0}^{t\wedge\theta_R}|f_1(s,x_s)||\bbE_s|ds  + 2M_t\\
&&+ C(1+|x_0|^{2l})\D^2 + \int_{0}^{t\wedge\theta_R} |g(\hat{s},s,y_{\hat{s}},y_s)-g(s,s,x_s,x_s)|^2 ds,
\eeao
where $M_t:=\int_{0}^{t\wedge\theta_R}|\bbE_s|\left( g(\hat{s},s,y_{\hat{s}},y_s)-g(s,s,x_s,x_s)\right)dW_s.$ It holds that
\beao
\bfE \sup_{0\leq t\leq T} |M_t| &\leq& 2\sqrt{32}\cdot\bfE \sqrt{\int_{0}^{T\wedge\theta_R}|\bbE_s|^2\left( g(\hat{s},s,y_{\hat{s}},y_s)-g(s,s,x_s,x_s)\right)^2ds}\\
&\leq&\bfE \sqrt{\sup_{0\leq s\leq T} |\bbE_{s\wedge\theta_R}|^2 \cdot128\int_{0}^{T\wedge\theta_R}\left( g(\hat{s},s,y_{\hat{s}},y_s)-g(s,s,x_s,x_s)\right)^2ds}\\
&\leq&\frac{1}{2}\bfE \sup_{0\leq s\leq T} |\bbE_{s\wedge\theta_R}|^2 + 64\bfE\int_{0}^{T\wedge\theta_R}\left( g(\hat{s},s,y_{\hat{s}},y_s)-g(s,s,x_s,x_s)\right)^2ds,
\eeao
thus we get that 
\beam
\nonumber\bfE\sup_{0\leq t \leq T}|\bbE_{t\wedge\theta_R}|^2 &\leq&2\bfE\sup_{0\leq t\leq T}\int_{0}^{t\wedge\theta_R}|f(\hat{s},s,y_{\hat{s}},y_s)-f(s,s,x_s,x_s)|^2ds +C_R\int_{0}^{t\wedge\theta_R} \bfE\sup_{0\leq l \leq s}|\bbE_l|^2 ds\\
\label{ExSD-eq:L2_conv_stoptime}&&+ C\D^2 + 130 \cdot \bfE\int_{0}^{T\wedge\theta_R} |g(\hat{s},s,y_{\hat{s}},y_s)-g(s,s,x_s,x_s)|^2ds.
\eeam
Assumption \ref{ExSD-ass:A} implies that 
$$
\int_{0}^{t\wedge\theta_R}|f_2(\hat{s},s,y_{\hat{s}},y_s)-f_2(s,s,x_s,x_s)|^2ds\leq \int_{0}^{t\wedge\theta_R}C_R\Big( |y_s - y_{\hat{s}}|^2 + |\bbE_s|^2  + |\hat{s}-s|^2\Big)ds  
$$
Moreover, it holds that
$$
\int_{0}^{t\wedge\theta_R}|\hat{s}-s|^2ds \leq \sum_{k=0}^{[t/\D-1]}\int_{t_k}^{t_{k+1}\wedge\theta_R} |t_k-s|^2ds.
$$
Taking the supremum over all $t\in[0,T]$ and then expectation we have
\beqq \label{ExSD-eq:L2_conv_stoptime_drift}
\bfE\sup_{0\leq t\leq T}\int_{0}^{t\wedge\theta_R}|f_2(\hat{s},s,y_{\hat{s}},y_s)-f_2(s,s,x_s,x_s)|^2ds\leq C_R\D + C_R\int_{0}^{T}\bfE\sup_{0\leq l\leq s}|\bbE_{l\wedge\theta_R}|^2ds  + C_R\D^2,
\eeqq
where we have used Lemma \ref{ExSD-lem:Error_SDbound} for $p=2.$ Using Assumption \ref{ExSD-ass:A} again we get that 
\beao
\!\!\!\!\!&&\!\!\!\!\!\!\int_{0}^{T\wedge\theta_R} \left( g(\hat{s},s,y_{\hat{s}},y_s)-g(s,s,x_s,x_s)\right)^2 ds \leq \int_{0}^{T\wedge\theta_R}C_R\Big( |y_s - y_{\hat{s}}|^2 + |\bbE_s|^2  + |\hat{s}-s|^2 + |y_{\hat{s}}-x_s|^{2q}\Big)ds\\
&\leq& \int_{0}^{T\wedge\theta_R}C_R\Big( |y_s - y_{\hat{s}}|^2 + |\bbE_s|^2  + |\hat{s}-s|^2 + \left(|y_{\hat{s}}-y_s|^{2q} + |\bbE_s|^{2q}\right)\Big)ds,
\eeao
where we have used the subadditivity property of $h(x)=x^{2q},$ thus taking expectation we have
\beao
&&\bfE\int_{0}^{T\wedge\theta_R} \left( g(\hat{s},s,y_{\hat{s}},y_s)-g(s,s,x_s,x_s)\right)^2 ds \leq C_R\D + C_R\int_{0}^{T}\bfE\sup_{0\leq l\leq s}|\bbE_{l\wedge\theta_R}|^2ds  + C_R\D^2\\
&&+ 2^{2q-1}C_RT\D^q + 2^{2q-1}C_R\int_{0}^{T\wedge\theta_R}(\bfE|\bbE_s|)^{2q}ds,
\eeao
where we have applied again Lemma \ref{ExSD-lem:Error_SDbound} for $p=2q$ and Jensen inequality with $2q<1.$ We get the following estimate
\beam\nonumber 
&&\bfE\int_{0}^{T\wedge\theta_R}(g(\hat{s},s,y_{\hat{s}},y_s)-g(s,s,x_s,x_s))^2ds \leq C_R\D^q  + C_R\int_{0}^{T}\bfE\sup_{0\leq l\leq s}(\bbE_{l\wedge\theta_R})^{2}ds\\ 
\nonumber&&+ C_R\int_{0}^{T}(K_{R,\D,m,q}(s))^{2q}ds\\
\label{ExSD-eq:L2_conv_stoptime_diff} &\leq& C_R\D^q  + C_R (K_{R,\D,m,q}(T))^{2q} + C_R\int_{0}^{T}\bfE\sup_{0\leq l\leq s}(\bbE_{l\wedge\theta_R})^{2}ds,
\eeam
where we have used Proposition \ref{ExSD-prop:L1_conv} and 
$$
K_{R,\D,m,q}(s):=\left( C_R\sqrt{\D} + \frac{C_R}{me_m}\D^{q} +\frac{C_R}{m}e_m +  e_{m-1}\right)e^{a_{R,m,q}s}.
$$
Plugging the estimates (\ref{ExSD-eq:L2_conv_stoptime_drift}), (\ref{ExSD-eq:L2_conv_stoptime_diff}) into (\ref{ExSD-eq:L2_conv_stoptime}) gives
\beao
\bfE\sup_{0\leq t \leq T}|\bbE_{t\wedge\theta_R}|^2 &\leq & C_R\D^q  + C_R (K_{R,\D,m,q}(T))^{2q} + C_R\int_{0}^{T}\bfE\sup_{0\leq l\leq s}(\bbE_{l\wedge\theta_R})^{2}ds\\
&\leq& \left(C_R\D^{q} + C_R (K_{R,\D,m,q}(T))^{2q}\right)e^{C_R T}\leq C_{R,\D,m}
\eeao
where we have applied the Gronwall inequality. Note that, given $R>0,$ the bound $C_{R,\D,m}$ can be made arbitrarily small by choosing big enough $m$ and small enough $\D.$ Relation (\ref{ExSD-eq:L2_conv}) becomes,
\beao
\bfE\sup_{0\leq t\leq T}|\bbE_t|^2 &\leq&C_{R,\D,m}  + 2^{\frac{3p-2}{p}}AR^{2-p}\\
&\leq& \un{C_{R,\D,m}}_{I_1}  + \un{2^{\frac{3p-2}{p}}AR^{2-p}}_{I_2}. 
\eeao
Given any $\vep>0,$ we may first choose $R$ such that  $I_2<\vep/2,$ then choose $m$ big enough and $\D$ small enough such that $I_1<\vep/2$ a concluding $\bfE\sup_{0\leq t\leq T}|\bbE_t|^2 <\vep$ as required to verify (\ref{ExSD-eq:strong_conv}).

\section{Proof of Theorem \ref{NSF-theorem:Pol rate}, Proposition \ref{NSF-prop:Pol rate} and Theorem \ref{NSF-theorem:Pol ratemulti}.}\label{NSF-sec:Proof}
\setcounter{equation}{0}

In this Section we prove our main strong convergence result. First, we provide uniform moment bounds for the original SDE and the SD scheme. We remind here that for notational reasons the processes $(W_t, x_t)$ stand for $(\wh{W}_t, \wh{x}_t).$  

\ble\label{NSF-lem:Moment Bounds}[Moment bounds for original problem and SD approximation]
Let Assumption \ref{NSF:assA} hold. Then 
$$
\bfE\sup_{0\leq t\leq T}|x_t|^p \bigvee \bfE\sup_{0\leq t\leq T}|y_t|^p\leq1,
$$
for any $p>0.$
\ele

\bpf[Proof of Lemma \ref{NSF-lem:Moment Bounds}]
The result is trivial since we already know that  $(x_t)$ satisfying (\ref{NSF-eq:WFmodel}) has the property $x_t\in D$ when $x_0\in D, D=(0,1),$ by Appendix \ref{NSF-ap:Boundary Classification} and regarding the bounds for the SD approximation it is clear, by its form (\ref{NSF-eq:SDsolution}), that they are valid. 
\epf

Now, let us rewrite the approximation process $(y_t^{SD})$
\beqq\label{NSF-eq:SD processNWagain}
y_t^{SD}= y_{t_n} + \left(k_1 - \frac{(k_3)^2}{4} + y_{t_n}(\frac{(k_3)^2}{2}-k_2)\right)\cdot\D  + \int_{t_n}^{t} \frac{(k_3)^2}{4}(1-2y_s)ds + k_3\int_{t_n}^{t} \sqrt{y_s(1-y_s)}dW_s.
\eeqq
In the general setting of (\ref{ExSD-eq:SD scheme}) we have 
$$
f_1(x)=k_1 - \frac{(k_3)^2}{4} + x\left(\frac{(k_3)^2}{2}-k_2\right), \qquad f_2(x)=\frac{(k_3)^2}{4}(1-2x), \qquad g(x)=b(x)=\sqrt{x(1-x)}.
$$
By the above representation, the form of the discretization  becomes apparent. We only discretized the drift coefficient of (\ref{NSF-eq:WFmodel}) in an additive way. Therefore, by an application of Theorem  \ref{ExSD-thm:strong_conv} we have the strong convergence result of Theorem \ref{NSF-theorem:Pol rate} 
$$
\lim_{\D\rightarrow0}\bfE\sup_{0\leq t\leq T}|y_t-x_t|^2=0.
$$

Now, we briefly sketch the proof of Proposition \ref{NSF-prop:Pol rate}. The process (\ref{NSF-eq:SD processNWalt}) is well-defined when $0<\wt{y}_{n}<1$ or equivalently  when $(k_3)^2<2k_2$ and $\D<-1/\beta$ using  (\ref{NSF-eq:initial_al}) and   (\ref{NSF-eq:additional_par}). The strong convergence result of Proposition \ref{NSF-prop:Pol rate} is a consequence of the triangle inequality and the following regularity-type result 
\beao
|\wt{y}_t^{SD} - y_t^{SD}| &=& \left|\sin^2 \left(\frac{k_3}{2}\D W_n^t + \arcsin(\sqrt{\wt{y}_n})\right) - \sin^2 \left(\frac{k_3}{2}\D W_n^t + \arcsin(\sqrt{y_n})\right)\right|\\
&\leq&2\left|\sin \left(\frac{k_3}{2}\D W_n^t + \arcsin(\sqrt{\wt{y}_n})\right) - \sin \left(\frac{k_3}{2}\D W_n^t + \arcsin(\sqrt{y_n})\right)\right|\\
&\leq&2\left| \arcsin(\sqrt{\wt{y}_n}) - \arcsin(\sqrt{y_n})\right|\\
&=&\left| \int_{\wt{y}_n}^{y_n}\frac{1}{\sqrt{z(1-z)}}dz \right|\leq |\alpha + \beta| \D \sup_{z\in \{ y_n, \wt{y}_n\}}\frac{1}{\sqrt{z(1-z)}}\leq C\cdot \D,
\eeao
for any $t\in (t_n, t_{n+1}],$ where $C$ is finite positive. 

Theorem \ref{NSF-theorem:Pol ratemulti} is an application of a slight generalization of Theorem  \ref{ExSD-thm:strong_conv} including multidimensional noise (see also \cite{halidias:2015}). Therefore, we omit the proof since one essentially repeats the steps of the proof of Theorem  \ref{ExSD-thm:strong_conv}. The auxiliary functions in the sense of  (\ref{ExSD-eq:SD scheme}) are, 
\beao
f_1(X)&=&k_1^{(1,1)}+k_1^{(1,2)}X^{(2)} - \frac{(k_3^{(1,1)})^2X^{(2)} + (k_3^{(1,2)})^2(1-X^{(1)}-X^{(2)})}{4(1-X^{(1)})} \\
&&+ X^{(1)}\left(\frac{(k_3^{(1,1)})^2X^{(2)} + (k_3^{(1,2)})^2(1-X^{(1)}-X^{(2)})}{2(1-X^{(1)})} -k_2^{(1)}\right)\\
f_2(X,Y^{(1)})&=&\frac{(k_3^{(1,1)})^2X^{(2)}+ (k_3^{(1,2)})^2(1-X^{(1)}-X^{(2)})}{4(1-X^{(1)})}(1-2Y^{(1)})\\
g_{11}(X,Y^{(1)})&=& k_3^{(1,1)}\sqrt{\frac{X^{(2)}}{1-X^{(1)}}}\sqrt{Y^{(1)}(1-Y^{(1)})}\\
g_{12}(X,Y^{(1)})&=& k_3^{(1,2)}\sqrt{\frac{1-X^{(1)}-X^{(2)}}{1-X^{(1)}}} \sqrt{Y^{(1)}(1-Y^{(1)})},
\eeao
for the evolution of the first component (\ref{NSF-eq:WFmultimodel}), where $X$ denotes the discretized part of the SDE, and accordingly for the second component (\ref{NSF-eq:WFmultimodel2})
\beao
f_1(X)&=&k_1^{(2,1)}+k_1^{(2,2)}X^{(1)} - \frac{(k_3^{(2,1)})^2X^{(1)} + (k_3^{(2,3)})^2(1-X^{(1)}-X^{(2)})}{4(1-X^{(2)})} \\
&&+ X^{(2)}\left(\frac{(k_3^{(2,1)})^2X^{(1)} + (k_3^{(2,3)})^2(1-X^{(1)}-X^{(2)})}{2(1-X^{(2)})} -k_2^{(2)}\right)\\
f_2(X,Y^{(2)})&=&\frac{(k_3^{(2,1)})^2X^{(1)}+ (k_3^{(2,3)})^2(1-X^{(1)}-X^{(2)})}{4(1-X^{(2)})}(1-2Y^{(2)})\\
g_{21}(X,Y^{(2)})&=& k_3^{(2,1)}\sqrt{\frac{X^{(1)}}{1-X^{(2)}}}\sqrt{Y^{(2)}(1-Y^{(2)})}\\
g_{23}(X,Y^{(2)})&=& k_3^{(2,3)}\sqrt{\frac{1-X^{(1)}-X^{(2)}}{1-X^{(2)}}}\sqrt{Y^{(2)}(1-Y^{(2)})},
\eeao
, see (\ref{NSF-eq:SDmultiprocess1}) and (\ref{NSF-eq:SDmultiprocess2}). By the above representation, the form of the discretization of (\ref{NSF-eq:WFmultimodel}) and (\ref{NSF-eq:WFmultimodel2}) becomes apparent. We discretized the drift coefficient in an additive and multiplicative way and the diffusion coefficient in a multiplicative way. 

\subsection*{Acknowledgements}
The author would like thank the anonymous referees for their helpful comments.

\bibliographystyle{unsrt}\baselineskip12pt 
\bibliography{wfmodel}

\appendix

\setcounter{argument}{0}
  
\section{Boundary classification of one-dimensional time-homogeneous SDEs.}\label{NSF-ap:Boundary Classification}
\setcounter{equation}{0}

Let us now recall some results \cite[Sec. 5.5]{karatzas_shreve:1988} concerning the boundary behavior of SDEs of the form,
\beqq\label{NSF-eq:General SDE in differential form}
 dX_t= a(X_t)dt + b(X_t)dW_t.
\eeqq
Let $I=(l,r)$ be an interval with $-\infty\leq l<r\leq \infty$ and define the exit time from $I$ to be 
$$ 
S:=\inf\{t\geq0: X_t\notin(l,r)\}.
$$
Let also the coefficients of (\ref{NSF-eq:General SDE in differential form}) satisfy the following conditions 
$$
b^2(x)>0, \quad \forall x\in I, \mbox{ \itshape(Non Degeneracy), (ND), }
$$
$$
\forall x\in I, \,\, \exists \ep>0: \int_{x-\ep}^{x+\ep}\frac{1 + |a(y)|}{b^2(y)}dy<\infty, \mbox{ \itshape(Local Integrability), (LI). }
$$
Then for $c\in I,$ we can define the scale function
\beqq\label{NSF-eq:scale function}
s(x):=\int_{c}^{x}e^{-2\int_c^y\frac{a(z)}{b^2(z)}dz}dy,
\eeqq
whose behavior at the endpoints of $I$ determines the boundary behavior of $(X_t)$  \cite[Prop. 5.5.22]{karatzas_shreve:1988}. In particular, we get that the dynamics (\ref{NSF-eq:WFmodel}) have a boundary behavior which is determined by the scale function
\beao
s(x)&=&\int_c^x \exp\big\{-2\int_c^y\frac{k_1-k_2z}{(k_3)^2z(1-z)}dz\big\}dy\\
&=&\int_c^x \exp\big\{-2\frac{k_1}{(k_3)^2}\int_c^y z^{-1}(1-z)^{-1}dz + 2\frac{k_2}{(k_3)^2}\int_c^y (1-z)^{-1}dz \big\}dy\\
&=&-\int_x^c\exp\big\{-2\frac{k_1}{(k_3)^2}\ln (y/c) + 2\frac{(k_1-k_2)}{(k_3)^2}\ln (1-y)/(1-c) \big\}dy \\
&=&-C\int_x^c y^{-2\frac{k_1}{(k_3)^2}}(1-y)^{2\frac{(k_1-k_2)}{(k_3)^2}}dy,
\eeao
for any $x\in I$  where $C>0.$  Let $I=(0,1)$ and take $c=1/2.$ We compute 
$$
s(0+)=-C\int_0^{1/2} y^{-2\frac{k_1}{(k_3)^2}}(1-y)^{2\frac{(k_1-k_2)}{(k_3)^2}}dy = -\infty,
$$
when $k_1>0$ and
\beao
s(1-) & = & C\int_{1/2}^1 y^{-2\frac{k_1}{(k_3)^2}}(1-y)^{2\frac{(k_1-k_2)}{(k_3)^2}}dy\\
& = & \infty,
\eeao
when $k_1<k_2,$ thus by \cite[Prop. 5.5.22a]{karatzas_shreve:1988} we have that $\bfP(S=\infty)=1$ that is $\bfP(0<x_t<1)=1.$

\section{Solution process of stochastic integral equations  (\ref{NSF-eq:SD process}), (\ref{NSF-eq:SDmultiprocess1}), (\ref{NSF-eq:SDmultiprocess2}).}\label{NSF-ap:Solution process}
\setcounter{equation}{0}

We will show that the process (\ref{NSF-eq:SDsolution}) for $n=0,$ is the solution of the stochastic integral equation  (\ref{NSF-eq:SD process}) for $n=0,$ that is 

\beqq\label{NSF-eq:SD solution n=0}
y_{t}^{SD}=\sin^2 \left( \frac{k_3}{2}W_t+ \arcsin(\sqrt{Y_0})\right),
\eeqq
satisfies
$$ 
y_t^{SD}= Y_0  + \int_{0}^{t} \frac{(k_3)^2}{4}(1-2y_s)ds + k_3\int_{0}^{t} \sqrt{y_s(1-y_s)}d\wh{W}_s,
$$ for $ t\in(0,t_{1}],$ with 
$$
Y_0:= x_0 + \left(k_1 - \frac{(k_3)^2}{4} + x_0\left(\frac{(k_3)^2}{2}-k_2\right)\right)\cdot\D\leq1.
$$
Relations  (\ref{NSF-eq:NewWiener}) and (\ref{NSF-eq:SDsgn term}) yield
 $$d\wh{W}_t:=\textup{sgn}(z_t)dW_t,$$
where  
$$
z_t= \sin \left(k_3\D W + 2\arcsin(\sqrt{Y_0})\right).
$$ 
The cases for $n=1,\ldots, N-1$ follow with the appropriate modifications.

We can write the increment of the Wiener process  as
$$
dW_t = 0\cdot dt + 1\cdot dW_t,
$$
and view (\ref{NSF-eq:SD solution n=0}) as a function of $W_t,$ i.e. $y=V(W)$ with 
\beao
\frac{dy}{dW}&=&2\sin \left(\frac{k_3}{2}\D W + \arcsin(\sqrt{Y_0})\right)\cos \left(\frac{k_3}{2}\D W + \arcsin(\sqrt{Y_0})\right)\cdot\frac{k_3}{2}\\
&=&\frac{k_3}{2}\sin \left(k_3\D W + 2\arcsin(\sqrt{Y_0})\right)\\
&=&\frac{k_3}{2}\sqrt{1-\cos^2 \left(k_3\D W + 2\arcsin(\sqrt{Y_0})\right)}\textup{sgn}\left[\sin \left(k_3\D W + 2\arcsin(\sqrt{Y_0})\right)\right]\\
&=&\frac{k_3}{2}\sqrt{1-(1-2y)^2}\textup{sgn}(z_t)\\
&=&k_3\sqrt{y(1-y)}\textup{sgn}(z_t)
\eeao
and 
\beao
\frac{d^2y}{dW^2}&=&\frac{(k_3)^2}{2}\cos \left(k_3\D W + 2\arcsin(\sqrt{Y_0})\right)\\
&=&\frac{(k_3)^2}{2}(1-2y).
\eeao
Application of It\^o's formula implies 
\beao
dy_t &=&  \frac{1}{2}V^{''}(W_t)dt + V^{'}(W_t)dW_t\\
&=& \frac{(k_3)^2}{4}(1-2y_t)  dt  + k_3\sqrt{y_t(1-y_t)}d\wh{W}_t.
\eeao
For the derivation of (\ref{NSF-eq:SDmultiprocess1}) and  (\ref{NSF-eq:SDmultiprocess2}) we now write the multidimensional Wiener process as 
$$
dW_t = \bbo_3\cdot dt + \bbI_3\cdot dW_t,
$$
where $\bbo_3$ is the zero $3\times3$ matrix and $\bbI_3$ the $3\times3$ identity matrix and apply appropriately the multidimensional It\^o formula.

\section{\texorpdfstring{Uniform moment estimate for $|y_t - x_t|^2$.}{Uniform moment estimate for abs(y(t) - x(t))2.}}\label{NSF-ap:uniform_moment_bound}
\setcounter{equation}{0}

In Theorem \ref{NSF-theorem:Pol rate} we actually proved that $\lim_{\D\downarrow0}\bfE\sup_{0\leq t\leq T}|y_{t}^{SD}-\wh{x}_{t}|^2=0.$ 
In order to finish the proof we have to find a uniform moment bound for $|x_t - \wh{x}_t|^2.$ In particular
applying  the triangle inequality 
\beao
\lim_{\D\downarrow0}\bfE\sup_{0\leq t\leq T}|y_{t}^{SD} - x_{t}|^2 &\leq&  \lim_{\D\downarrow0}\bfE\sup_{0\leq t\leq T}|y_{t}^{SD}-\wh{x}_{t}|^2 + \lim_{\D\downarrow0}\bfE\sup_{0\leq t\leq T}|x_t-\wh{x}_{t}|^2\\
&\leq&\lim_{\D\downarrow0}\bfE\sup_{0\leq t\leq T}|x_t-\wh{x}_{t}|^2\eeao
thus it suffices to show
$$
\lim_{\D\downarrow0}\bfE\sup_{0\leq t\leq T}|x_t-\wh{x}_{t}|^2=0; 
$$
this follows by (\ref{NSF-eq:WFmodel}) and (\ref{NSF-eq:WFmodel_NW}).

\end{document}